\documentclass[sort&compress,preprint]{elsarticle}

\usepackage[utf8]{inputenc} 
\usepackage{amssymb}
\usepackage{amsfonts}
\usepackage{mathtools}
\usepackage{bbm}
\usepackage{float}
\usepackage{multicol}
\usepackage{placeins}
\usepackage[caption=false]{subfig}
\usepackage{graphicx}
\usepackage[inline]{enumitem}
\usepackage{lineno}
\usepackage{hyperref}
\usepackage{cleveref}
\usepackage{algorithmic}
\usepackage{float}
\usepackage{csl}

\newcommand{\DKL}{\operatorname{D}_{\rm KL}}
\newcommand{\DhKL}{\widehat{\operatorname{D}}_{\rm KL}}

\newif\ifcomplete
\newif\ifdraft
\newif\ifreport

\newcommand{\sfrac}[2]{\mbox{\footnotesize$\displaystyle\frac{#1}{#2}$}} 
\newcommand{\nablax}{\nabla_{\mkern-4mu\x}}
\newcommand{\nablaz}{\nabla_{\mkern-4mu\z}}

\newcommand{\Hobs}{\mathcal{H}}
\newcommand{\HH}{\mathbf H}
\newcommand{\Nobs}{ {{\rm N}_{\rm obs}} }

\newcommand{\C}{\mathbf{C}}

\newcommand{\Id}{ \mathbf{I} }

\renewcommand{\P}{\mathbf{P}}

\newcommand{\Pb}{\mathbf{P}^{\rm b}}

\newcommand{\Model}{\mathcal{M}}

\newcommand{\Mtlm}{\mathbf{M}}
\newcommand{\Madj}{\mathbf{M}^*}

\newcommand{\tr}{^{\mkern-1.5mu\mathrm{T}}} 

\newcommand{\x}[1][]{%
   \ifthenelse{ \equal{#1}{} }
      {\mathbf{x}}
      {\mathbf{x}^{[#1]}}}
\newcommand{\xb}[1][]{%
   \ifthenelse{ \equal{#1}{} }
      {\mathbf{x}^{\rm b}}
      {\mathbf{x}^{{\rm b}[#1]}}}
\newcommand{\xa}[1][]{%
   \ifthenelse{ \equal{#1}{} }
      {\mathbf{x}^{\rm a}}
      {\mathbf{x}^{{\rm a}[#1]}}}
\newcommand{\xf}[1][]{%
   \ifthenelse{ \equal{#1}{} }
      {\mathbf{x}^{{\rm f}}}
      {\mathbf{x}^{{\rm f}[#1]}}}
\newcommand{\y}[1][]{%
   \ifthenelse{ \equal{#1}{} }
      {\mathbf{y}}
      {\mathbf{y}^{[#1]}}}
\newcommand{\xprop}[1][]{
   \ifthenelse{ \equal{#1}{} }
      {\mathbf{x}^{\rm q}}
      {\mathbf{x}^{{\rm q}[#1]}}}

\newcommand{\bigdot}[1]{\accentset{\mbox{\large\bfseries .}}{#1}}
\newcommand{\dx}[1][]{%
   \ifthenelse{ \equal{#1}{} }
      {\bigdot{\mathbf{x}}}
      {\bigdot{\mathbf{x}}^{[#1]}}}
\newcommand{\dxb}[1][]{%
   \ifthenelse{ \equal{#1}{} }
      {\bigdot{\mathbf{x}}^{\rm b}}
      {\bigdot{\mathbf{x}}^{{\rm b}[#1]}}}
\newcommand{\dxa}[1][]{%
   \ifthenelse{ \equal{#1}{} }
      {\bigdot{\mathbf{x}}^{\rm a}}
      {\bigdot{\mathbf{x}}^{{\rm a}[#1]}}}

\newcommand{\xt}{ \mathbf{x}^{\rm true} }

\newcommand{\xmean}{ \overline{\mathbf{x}} }
\newcommand{\xmeanb}{ \overline{\mathbf{x}}^{\rm b} }
\newcommand{\xmeanf}{ \overline{\mathbf{x}}^{\rm f} }
\newcommand{\xmeana}{ \overline{\mathbf{x}}^{\rm a} }

\newcommand{\hofx}[1][]{%
   \ifthenelse{ \equal{#1}{} }
      {\mathbf{z}}
      {\mathbf{z}^{[#1]}}}
\newcommand{\hofxb}[1][]{%
   \ifthenelse{ \equal{#1}{} }
      {\mathbf{z}^{\rm b}}
      {\mathbf{z}^{{\rm b}[#1]}}}
\newcommand{\hofxa}[1][]{%
   \ifthenelse{ \equal{#1}{} }
      {\mathbf{z}^{\rm a}}
      {\mathbf{z}^{{\rm a}[#1]}}}
\newcommand{\dhofxb}[1][]{%
   \ifthenelse{ \equal{#1}{} }
      {\bigdot{\mathbf{z}}^{\rm b}}
      {\bigdot{\mathbf{z}}^{{\rm b}[#1]}}}
\newcommand{\dhofxa}[1][]{%
   \ifthenelse{ \equal{#1}{} }
      {\bigdot{\mathbf{z}}^{\rm a}}
      {\bigdot{\mathbf{z}}^{{\rm a}[#1]}}}
\renewcommand{\d}[1][]{%
   \ifthenelse{ \equal{#1}{} }
      {\mathbf{d}}
      {\mathbf{d}^{[#1]}}}

\newcommand{\z}{ \mathbf{z} }

\newcommand{\errb}[1][]{%
   \ifthenelse{ \equal{#1}{} }
      {\varepsilon^{\rm b}}
      {\varepsilon^{{\rm b}[#1]}}}
\newcommand{\erra}[1][]{%
   \ifthenelse{ \equal{#1}{} }
      {\varepsilon^{\rm a}}
      {\varepsilon^{{\rm a}[#1]}}}
\newcommand{\erro}[1][]{%
   \ifthenelse{ \equal{#1}{} }
      {\varepsilon^{\rm obs}}
      {\varepsilon^{{\rm obs}[#1]}}}

\newcommand{\errm}[1][]{%
   \ifthenelse{ \equal{#1}{} }
      {\eta}
      {\eta^{[#1]}}}
\newcommand{\wb}[1][]{%
   \ifthenelse{ \equal{#1}{} }
      {w^{\rm b}}
      {w^{{\rm b}[#1]}}}
\newcommand{\wa}[1][]{%
   \ifthenelse{ \equal{#1}{} }
      {w^{\rm a}}
      {w^{{\rm a}[#1]}}}
\newcommand{\wprop}[1][]{
   \ifthenelse{ \equal{#1}{} }
      {w^{\rm a/\!q}}
      {w^{{\rm a/\!q}[#1]}}}

\renewcommand{\Re}{\mathbbm{R}}

\newcommand{\Prob}{\mathcal{P}}
\newcommand{\PA}{ {\mathcal P}^{\rm a} }
\newcommand{\PB}{ {\mathcal P}^{\rm b} }

\newcommand{\PO}{ {\mathcal P}^{\rm obs} }

\newcommand{\expectw}[2]{\mathtt{E}_{#1}\left[#2\right]}

\newcommand{\cov}[1]{\mathtt{Cov}\left[#1\right]}

\newcommand{\from}{\sim}

\renewcommand{\mid}{\,|\,}

\newcommand{\norm}[1]{\bigl\Vert #1 \bigr\Vert}

\newcommand{\dif}{\mathrm{d}}

\newcommand{\Xens}{\mathbf{X}}

\newcommand{\Xensb}{\mathbf{X}^{\rm b}}
\newcommand{\Xensa}{\mathbf{X}^{\rm a}}

\newcommand{\Nens}{{\rm N}_{\rm ens}}
\newcommand{\Nstate}{{\rm N}_{\rm state}}

\newtheorem{thm}{Theorem}

\newdefinition{rmk}{Remark}
\newproof{pf}{Proof}
\newproof{pot}{Proof of Theorem \ref{thm2}}
\newtheorem{assumption}{Assumption}
\newtheorem{example}{Example}
\newtheorem{remark}{Remark}

\crefname{thm}{Theorem}{Theorems}
\crefname{assumption}{Assumption}{Assumptions}

\newcommand{\titlestring}{Ensemble Variational Fokker-Planck Methods for Data Assimilation}

\newcommand{\authorstring}{Amit N Subrahmanya, Andrey A Popov, Adrian Sandu}
\newcommand{\emailstring}{amitns@vt.edu, apopov@vt.edu, sandu@vt.edu}

\def\*#1{\boldsymbol{\mathbf{#1}}}
\def\!#1{\mathcal{#1}}

\newcommand{\sig}{\boldsymbol{\sigma}}
\newcommand{\syt}{\tau}
\newcommand{\syth}{\zeta}
\newcommand{\diffu}{\alpha}
\newcommand{\reg}{\beta}
\newcommand{\dive}{\operatorname{div}}

\newcommand{\xsyt}[1][]{
   \ifthenelse{ \equal{#1}{} }
      {\mathbf{x}_{\rm \syt}}
      {\mathbf{x}_{{\rm \syt}[#1]}}}
  
\newcommand{\Xenssyt}[1][]{
   \ifthenelse{ \equal{#1}{} }
      {\mathbf{X}_{\rm \syt}}
      {\mathbf{X}_{{\rm \syt}[#1]}}}

\reporttrue

\journal{Journal of Computational Physics}

\begin{document}

\csltitle{\titlestring}
\cslauthor{\authorstring}
\cslyear{21}
\cslreportnumber{10}
\cslemail{\emailstring}
\csltitlepage

\begin{frontmatter}

\title{\titlestring}

\author[1]{Amit N. Subrahmanya\corref{cor1}}
\ead{amitns@vt.edu}
\author[1]{Andrey A. Popov}
\ead{apopov.vt.edu}
\author[1]{Adrian Sandu}
\ead{asandu7@vt.edu}

\cortext[cor1]{Corresponding author}

\affiliation[1]{organization={Computational Science Laboratory, Department of Computer Science, Virginia Tech},
addressline={620 Drillfield Dr.},
city={Blacksburg},
postcode={24061},
state={Virginia},
country={USA}}

\begin{abstract}
Particle flow filters solve Bayesian inference problems by smoothly transforming a set of particles into samples from the posterior distribution. Particles move in state space under the flow of an McKean-Vlasov-It\^{o} process.
This work introduces the Variational Fokker-Planck (VFP) framework for data assimilation, a general approach that includes previously known particle flow filters as special cases.
The McKean-Vlasov-It\^{o} process that transforms particles is defined via an optimal drift that depends on the selected diffusion term. 
It is established that the underlying probability density - sampled by the ensemble of particles - converges to the Bayesian posterior probability density.  
For a finite number of particles the optimal drift {contains} a regularization term that nudges particles toward becoming independent random variables. 
Based on this analysis, we derive computationally-feasible approximate regularization approaches that penalize the mutual information between pairs of particles, and avoid particle collapse. 
Moreover, the diffusion plays a role akin to a particle rejuvenation approach that aims to alleviate particle collapse.  
The VFP framework is very flexible. 
Different assumptions on prior and intermediate probability distributions can be used to implement the optimal drift, and localization and covariance shrinkage can be applied to alleviate the curse of dimensionality. 
A robust implicit-explicit method is discussed for the efficient integration of stiff McKean-Vlasov-It\^{o} processes.
The effectiveness of the VFP framework is demonstrated on three progressively more challenging test problems, namely the Lorenz\,'63, Lorenz\,'96 and the quasi-geostrophic equations.
\end{abstract}

\ifreport

\fi

\begin{keyword}
Bayesian Inference \sep Data Assimilation \sep Particle Filters \sep Particle Flow
\MSC 65C05 \sep 93E11 \sep 62F15 \sep 86A22

\end{keyword}

\end{frontmatter}


\section{Introduction}
\label{sec:intro}

Data assimilation (DA)~\cite{Asch_2016_book,Reich_2015_book} seeks to estimate the state of a physical system by optimally combining sparse and noisy observations of reality with background information obtained from a computational model of the system.
As exact Bayesian inference for this state estimation problem is computationally intractable, statistical sampling methods are frequently used to perform approximate inference. 

State-of-the-art statistical methods for data assimilation include the Ensemble Kalman Filter (EnKF)~\cite{Evensen_1994, Evensen_2003, Burgers_1998_EnKF}, and its variants such as the Ensemble Transform Kalman Filter~\cite{Bishop_2001_ETKF}, the Ensemble Kalman Smoother (EnKS)~\cite{Evensen_2000}, all of which make Gaussian assumptions on the distribution of the samples. 
One issue with the aforementioned methods is particle collapse -- when the samples in an ensemble become similar to each other, the ensemble covariance becomes small, and the filter trusts the model while discarding the observations, which leads to the divergence of the analysis trajectory from the truth.
Popular heuristics to prevent divergence include covariance inflation~\cite{Anderson_1999_MC-implementation}, and particle rejuvenation~\cite{Reich_2016_hybrid,popovamit}.
Another issue with these methods is the curse of dimensionality~\cite{Hastie_2001_statsbook} -- statistical estimates of the covariance are low rank and inaccurate due to a dearth of samples in high dimensions, requiring heuristic corrections such as covariance localization~\cite{Anderson_2007_localization} and covariance shrinkage~\cite{Chen_2010_shrinkage}. 

Particle filters~\cite{vanLeeuwen_2009_PF-review,vanLeeuwen_2019_PF-review} make little to no assumptions about any of the underlying distributions; distributions are represented empirically by {an} ensemble of particles, each with a certain weight quantifying its likelihood.
The inference step updates the weights rather than updating the particle states.
Particle filters suffer from degeneracy -- when the weights of a small subset of particles are large, and the weights of the remaining particles are close to zero, the effective number of samples is small and the accuracy of the filter is degraded.
Particles must be resampled periodically to avoid degeneracy.
While attractive due to their generality, traditional particle filters are impractical for usage in high dimensional problems as they require exceedingly large numbers of particles.
More robust approaches to particle filtering have been recently developed based on probability transport maps, and include the Ensemble Transport Particle Filter (ETPF)~\cite{Reich_2013_ETPF}, the second-order Ensemble Transport Particle Filter(ETPF2)~\cite{Reich_2017_ETPF_SOA}, a coupling-based ensemble filter~\cite{Marzouk_2019_nonlinear-coupling} and the Marginal Adjusted Rank Histogram Filter (MARHF)~\cite{Anderson_2020_rank-filter}.

Particle flow filtering, where particles move continuously in the state space toward samples from a posterior distribution, has attracted considerable attention recently as a general methodology for Bayesian inference. 
The particle motion is governed by the flow of a differential equation. To define this flow, the Stein variational gradient descent method~\cite{Liu_2016_Stein} uses the equality between the Stein discrepancy and the gradient of the Kullback-Leibler (KL) divergence~\cite{Kolmogoroff_1931_FPE} between the current and posterior distributions. The aim is to progressively minimize the KL divergence between the posterior distribution and the sequence of intermediate particle distributions and the posterior distribution.
A closed form solution to the flow is defined by embedding the particles into a reproducing kernel Hilbert space (RKHS), whose kernel must be meticulously chosen.
The mapping particle filter (MPF)~\cite{Pulido_2019_mapping-PF} employs the Stein variational gradient descent approach to perform data assimilation. 
{Scalability of MPF to higher dimensions is challenging, as MPF is biased for a small ensemble spread due to a finite number of samples}.
{The particle flow filter (PFF)~\cite{Hu_2020_mapping-PF} reduces this bias by employing different localized kernel functions for each state variable, but the problem does not disappear.}

While the previously discussed approaches~\cite{Liu_2016_Stein,Pulido_2019_mapping-PF,Hu_2020_mapping-PF} view the particle flow problem through the lens of optimization, other works~\cite{Reich_2019_discrete-gradients,Reich_2021_FokkerPlanck} take a dynamical system point of view where a McKean-Vlasov-It\^{o} process evolves particles {from sampling a prior distribution toward sampling a target (or posterior) distribution}.
Reich and Weissmann~\cite{Reich_2021_FokkerPlanck} consider the dynamics of an interacting system of particles, and the evolution of the corresponding probability distributions via the Fokker Planck equation. 
They discuss sufficient conditions that lead to the convergence of the evolving distribution of samples to the posterior distribution, which allows to perform Bayesian inference with a wide variety of interacting particle approximations.
A major drawback of this approach is that the evolution of the interactive particle system is highly stiff, and requires expensive numerical integration approaches. 
A related approach proposed by Garbuno-Inigo et. al.~\cite{Stuart_2020_gradient-EnKF} uses interacting Langevin diffusions to define particle flows.
The Fokker-Planck equation associated with the stochastic process has an exploitable gradient structure built on the Wasserstein metric and the covariance of the diffusion, which ensures convergence to the desired posterior distribution.
A derivative-free implementation of the dynamics is proposed, which allows to extend Ensemble Kalman Inversion~\cite{Stuart_2013_EnKF_inversion} to compute samples from a Bayesian posterior.  

This work introduces a generalized variational Fokker-Planck (VFP) approach to data assimilation.
Much of the theory for stochastic processes moving a probability density to a desired target via Fokker-Planck dynamics has been discussed by Jordan et al ~\cite{Jordan_1998_FokkerPlanck}.
We show that previously described methods such as the MPF, the PFF, and {(first-order, overdamped)} Langevin-based filters are in fact, particular formulations in the VFP framework.
{Specifically, the deterministic formulations of MPF and PFF can be obtained by embedding the particle dynamics in a reproducing kernel Hilbert space(RKHS), with no diffusion.} 
The VFP framework also extends Fokker-Planck~\cite{Reich_2021_FokkerPlanck} and Langevin dynamics~\cite{Stuart_2020_gradient-EnKF} filters, and offers wider flexibility in defining particle flows.

The key contributions of this paper are: 
\begin{enumerate*}[label={(\roman*)}]
	\item a generalized formulation of the variational Fokker-Planck framework that subsumes multiple previously proposed ensemble variational  data assimilation methods,
	\item derivation of the optimal drift of a McKean-Vlasov-It\^{o} process -- that depends on the selected diffusion term -- to push particles towards the posterior distributions,
    \item a general implementation of VFP via combinations of parameterized distributions, 
	\item derivation of regularization terms to ensure particle diversity by nudging particles toward becoming independent random variables,
    \item an extension of the formalism to solve smoothing problems in both strong-constraint (perfect model) and for weak-constraint (model with errors) cases,
    \item inclusion of localization and covariance shrinkage in the VFP approach for high dimensional problems,
	\item discussion of a partitioned linearly-implicit-explicit stochastic time integration method to evolve stiff McKean-Vlasov-It\^{o} processes.
\end{enumerate*}

The remainder of this paper is organized as follows. The general discrete time data assimilation problem is reviewed in \Cref{sec:background} along with a description of notation.
\Cref{sec:vfpd} develops the proposed variational Fokker-Planck framework, including the derivation of the optimal drift and regularization terms, options to parametrize the intermediate distributions, and implementation aspects.
Examples of particular VFP filters and smoothers are shown in \Cref{sec:examples-vfp}.
\ifreport
The importance of regularization and diffusion is illustrated with the help of an example in \Cref{sec:illustrate-regularization-diffusion}.\fi
Application of localization and covariance shrinkage in the VFP framework are discussed in \Cref{sec:apply-localization}.
Numerical experiments to validate the methodology using Lorenz\,'63~\cite{Tandeo_2015_l63,Lorenz_1963_L63}, Lorenz\,'96~\cite{Lorenz_1996_L96,vanKekem_2018_l96dynamics} and quasi-geostrophic equations~\cite{Charney_1947_QG, San_2015_qge} are reported in \Cref{sec:expts}.
Concluding remarks are drawn in \Cref{sec:conc}.

\section{Background}
\label{sec:background}

Science seeks to simulate and forecast natural processes that evolve in time, with dynamics that are too complex to be fully known.
Let $\xt_k \in \Re^{\Nstate}$ denote the true state of a dynamical system at time $t_k$, representing (a finite dimensional projection of) the state of the natural process.
Due to our lack of knowledge, our simulation represents only an estimate of the truth.
The background (prior) estimate of the state is represented by a random variable $\xb_k \sim \PB_k$ whose distribution quantifies the prior uncertainty in our knowledge.
Here, $\Prob$ represents a probability density. 
State estimates are evolved in time by the computational model $\Model_{k,k+1} : \Re^{\Nstate} \to \Re^{\Nstate}$ which does not  fully capture the dynamics of the natural process, therefore the simulation results will slowly diverge from reality. 
To prevent this, our computed estimate must be combined with observations of the true state. 
An observation is defined as 
\begin{equation}\label{eq:observation}
    \y_k = \Hobs_k (\xt_k) + \erro_k, \quad \y_k \in \Re^{\Nstate}, \quad k \ge 0,
\end{equation}
where $\Hobs_k: \Re^{\Nstate} \to \Re^{\Nobs}$ is a non-linear observation operator.
It is assumed that the observation is corrupted by observation errors from a known error distribution $\erro_k \from \PO_k$, and that observation errors at different times are {independent random} variables.
In most operational problems, observations are spatially sparse i.e. $\Nobs \ll \Nstate$ as they are expensive to acquire.

Our goal is to perform Bayesian inference using these two sources of information  --  the background and observation  --  to decrease uncertainty and obtain an improved estimate of the true state.
This improved estimate is another random variable $\xa_k \sim \PA_k$, called the analysis (posterior), that represents our total knowledge about the state at time $t_k$. The posterior distribution given by Bayes' rule \cite{Robert_2004_Monte} is 
\begin{equation}\label{eq:Bayes-rule}
    \PA_k(\x_k) = \PB_k(\x_k | \y_k) = \frac{\PO_k(\y_k | \x_k) \, \PB_k(\x_k)}{\Prob_k(\y_k)}.
\end{equation}
By propagating the analysis in time from $t_k$ to $t_{k+1}$ through the model operator:
\begin{equation}\label{eq:model}
    \xb_{k+1} = \Model_{k,k+1}(\xa_k) + {\errm_{k + 1}}, \quad k \ge 0,
\end{equation}
a new prior is obtained at time $t_{k+1}$, and the cycle can begin anew. 
{In this paper the model error term $\errm_{k+1} \sim \Prob_{k+1}^{\Model}(\errm)$ is taken to be zero, i.e., we assume a perfect model (also referred to as a strong constraint in certain variational applications\cite{Evensen_2022_book,Sandu_2011_assimilationOverview})}.

Consider the data assimilation window $[t_0,t_{K}]$.
The filtering approach to the Bayesian inference problem~\cref{eq:Bayes-rule} incorporates only the observation at the current time, and sequentially produces analyses $\xa_k$ conditioned by all past observations $\y_{0:k}$ for $0 \le k \le K$: 
\begin{equation}\label{eq:Bayes-filter}
   \xa_k \from \Prob(\x_k \mid \y_{0:k}) = \frac{\PO_{k}(\y_{k} \mid \x_k)}{ \Prob_k(\y_{k})}\, \Prob_k(\x_k \mid \y_{0:k-1})
   \stackrel{\errm_i=0}{=} \left[\prod_{i = 0}^k \frac{\PO_{i}(\y_{i} \mid \x_i)}{ \Prob_i(\y_{i})}\right]\PB_0(\x_0).
\end{equation}
In contrast, the {strong constraint} smoothing approach incorporates all present and future observations within the assimilation window $[t_0,t_{K}]$ into the current analysis starting from $t_0$ as, 
\begin{equation}\label{eq:Bayes-smooth}
   \xa_0 \from \Prob(\x_0 \mid \y_{0:K}) 
   \stackrel{\errm_i=0}{=} \left[\prod_{i = 0}^K \frac{\PO_{i}(\y_{i} \mid \x_i)}{ \Prob_i(\y_{i})}\right]\PB_0(\x_0)
\end{equation}
{Additionally, the weak constraint~\cite{Evensen_2022_book} smoothing approach considers the (non-zero) model error distributions over the assimilation window as, 
\begin{equation}\label{eq:weak-Bayes-smooth}
   \xa_0 \from \Prob(\x_0 \mid \y_{0:K}) =
   \left[\prod_{i = 1}^K \frac{\PO_{i}(\y_{i} \mid \x_i) \Prob_{i}^{\Model}(\errm_i) }{ \Prob_i(\y_{i})}\right] \left( \frac{\PO_{0}(\y_{0} \mid \x_0)}{ \Prob_0(\y_{0})} \PB_0(\x_0) \right)
\end{equation}}
In practice, performing exact Bayesian inference {as in} \cref{eq:Bayes-filter,eq:Bayes-smooth,eq:weak-Bayes-smooth} is computationally infeasible. 
Most tractable methods work via Monte-Carlo approaches that represent the probability densities used in inference \cref{eq:Bayes-rule} by their empirical counterparts.
To this end, we denote an ensemble of $\Nens$ realizations (or samples) of the random state variable $\x \from \Prob$ as
\begin{equation}\label{eq:ensemble-definition}
    \Xens \coloneqq \left[ \x[1], \x[2], \cdots, \x[\Nens] \right] \in \Re^{\Nstate \times \Nens}.  
\end{equation} 
In the ensemble limit of $\Nens \to \infty$, the empirical measure distribution of the ensemble,
\begin{equation}
    \widetilde{\Prob}(\x) = \frac{1}{\Nens}\sum_{i=1}^{\Nens} \delta_{\x^{[i]}}(\x),
\end{equation}
converges weakly, almost surely to the distribution of the random variable $\Prob(\x)$.
When the random variable $\x$ describes the state of a dynamical system, each ensemble member $\x[e]$ is also called a particle to hint at its propagation in time. 
Particles are used to estimate statistics of the probability distribution $\Prob$. For example, the empirical mean $\xmean$, anomalies $\*A$, and covariance $\P$ are defined as
\begin{equation}\label{eq:ensemble-mean-and-covariance}
    \xmean = \frac{1}{\Nens}\,\Xens\*1_{\Nens}, \quad \*A = \frac{1}{\sqrt{\Nens - 1}} (\Xens - \xmean\,\*1_{\Nens}\tr),\quad \P = \*A \,{\*A}\tr,
\end{equation}
respectively, where $\*1_p$ represents a $p$-dimensional vector of ones.
The background ensemble of particles $\Xensb_k = \left[ \xb[1]_k,  \cdots, \xb[\Nens]_k \right]$ represents the background probability density as $\xb[e]_k \from \PB_k(\x)$.
The data assimilation problem now is to produce an analysis ensemble $\Xensa_k = \left[ \xa[1]_k,  \cdots, \xa[\Nens]_k \right]$ with $\xa[e]_k \from \PA_k(\x)$ that represents the analysis probability density.
For the remainder of this paper we ignore the physical time subscripts: $\xb \equiv \xb_k$, $\PB \equiv \PB_k$, $\xa \equiv \xa_k$, $\PA \equiv \PA_k$, unless necessary.
%

\section{The Variational Fokker-Planck approach to data assimilation}
\label{sec:vfpd}

Variational particle filters use a dynamical system approach to transform a set of Monte Carlo samples $\xb[e]$ from the prior distribution into samples $\xa[e]$ from the posterior distribution \cite{Taghvaei_2019_particle-flow-accelerated,Jordan_1998_FokkerPlanck,Stuart_2020_gradient-EnKF,Daum_2011_particle-flow,Reich_2021_FokkerPlanck,Reich_2010_localization}.
Here, the particles move in state space, according to a differential equation in artificial time $\syt$, such that the underlying probability distributions evolve from the prior to the posterior.
Two approaches to moving the particles have been proposed. One approach formulates a flow over a finite time interval $\syt \in [0,1]$, starts with the prior distribution at $\syt = 0$, and reaches the posterior distribution at $\syt = 1$ \cite{Reich_2012_Gauss-mixture}. 
A second approach formulates a flow over an interval $\syt \in [0,\infty)$, and maps any initial distribution to the posterior distribution asymptotically when $\syt \to \infty$. 
Examples of such filters include the Stein variational gradient descent \cite{Liu_2016_Stein,Liu_2017_Stein}, the mapping particle filter \cite{Hu_2020_mapping-PF,Pulido_2019_mapping-PF}, interacting Langevin diffusions \cite{Stuart_2020_gradient-EnKF}, and Fokker-Planck particle systems \cite{Reich_2021_FokkerPlanck}.  
In this work, we generalize the second approach and propose the Variational Fokker-Planck (VFP) framework for data assimilation. 

\begin{figure}[tbhp]
    \centering
    \includegraphics[width=0.4\linewidth]{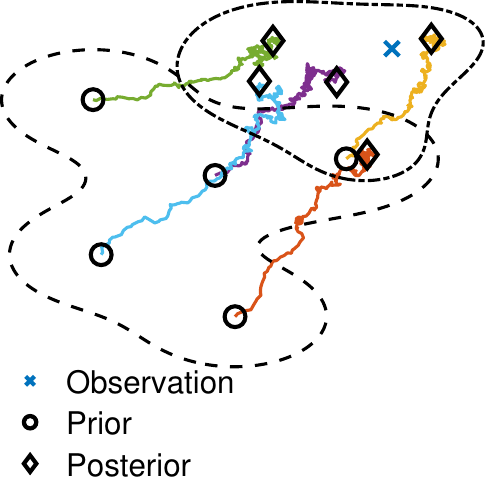}
    \caption{Particles sampled from a prior distribution move toward samples from the posterior distribution under the flow of a stochastic differential equation.}
    \label{fig:density-flows}
\end{figure}

The main idea of the VFP approach is as follows.
The initial configuration of the system is a set of particles \cref{eq:ensemble-definition} drawn from the prior distribution $\xb[e] \sim \PB$.
The particles move under the flow of a McKean-Vlasov-It\^{o} process, and their underlying distribution evolves according to the corresponding Fokker-Planck equation~\cite{Kolmogoroff_1931_FPE}. 
This McKean-Vlasov-It\^{o} process is defined such as to push the particles towards becoming samples of the posterior distribution $\xa[e] \sim\PA$. 
{This ensemble $\Xenssyt$, evolving in synthetic time $\syt$, is referred to as the current (or intermediate) ensemble; each particle from this ensemble is a sample of the current (or intermediate) distribution, i.e. $\xsyt^{[e]} \sim q_\syt$.}
This idea is illustrated in \cref{fig:density-flows}. 
The prior/background particles are depicted by the circles whose density is shown by the dashed line. 
These particles flow towards the observation under an optimal McKean-Vlasov-It\^{o} process depicted by colored lines. 
The final positions of the particles are marked by diamonds, with the dash-dotted line around the analysis particles representing the posterior distribution. 
%

\subsection{Derivation of the optimal drift}  
\label{subsec:VFP-derivation}

{We focus on the case where the posterior probability is absolutely continuous with {respect to} the Lebesgue measure, and the following assumption hold.
\begin{assumption}
\label{ass:connectivity}
{Let $\Omega = \Re^{\Nstate}$.} Consider the set of smooth probability densities \cite{Stuart_2020_gradient-EnKF} :
\begin{equation}
\label{eq:smooth-pdfs}
\begin{split}
\Pi &= \Big\{ q \in C^{\infty}(\Omega)
~:~ q(\x) > 0~a.e., ~ \int_\Omega q(\x) d\x = 1, ~  \int_\Omega \Vert \x \Vert^2 q(\x) d\x < \infty \Big\}.
\end{split}
\end{equation}
We make the following assumptions.
\begin{itemize}
\item $q_0,\PA \in \Pi$.
\item The solutions of the Fokker-Planck equation, as in \cref{eq:fokker-planck-vlasov}, are smooth, $q_\tau \in C^1\big( [0,\infty), \Pi \big)$. 
\end{itemize}
\end{assumption}
}
Consider an initial value McKean-Vlasov-It\^{o} process~\cite{Kloeden_2011_sdebook,Evans_2012_sdebook,Barbu_2010_FokkerPlanck} that acts on a random variable $\xsyt \in \Re^{\Nstate}$ evolving in artificial-time $\syt$:
\begin{equation}\label{eq:ito-process-difussion}
    \dif \xsyt  = \mathbf{F}(\syt, \xsyt, q_\syt)\, \dif \syt + \sig(\syt, \xsyt, q_\syt)\, \dif \mathbf{W}_\syt, \quad \syt \geq 0, \quad \x_0 \from q_0,
\end{equation}
where $\mathbf{F}: \Re_+ \times \Re^{\Nstate} \times \Pi \to \Re^{\Nstate}$ is a smooth drift term, $\sig : \Re_+ \times  \Re^{\Nstate}  \times \Pi \to \Re^{\Nstate \times M}$ is a smooth diffusion matrix, $q_\syt \in \Pi$, and $\mathbf{W}_\syt \in \Re^{M}$ is an $M$-dimensional standard Wiener process. 
We make the assumption that $\sig(t,\x,q)$ is a functional of $q$, i.e., $\sig$ depends on the state $\x$ only via the second argument.
The random variable evolved by \cref{eq:ito-process-difussion} in artificial-time has a probability density $\xsyt \from q_\syt$, which is the solution of the corresponding Fokker-Planck-Vlasov equation \cite{Barbu_2010_FokkerPlanck}:
\ifreport
\begin{equation}
\label{eq:fokker-planck-vlasov2}
\begin{split}
    \frac{\partial q_\syt(\x)}{\partial \syt} &= -\sum_{i = 1}^{\Nstate} \frac{\partial}{\partial \x_i} \big( q_\syt(\x) \mathbf{F}_i(\syt, \x, q_\syt) \big) \\
    &\quad + \sum_{i = 1}^{\Nstate} \sum_{j = 1}^{\Nstate} \frac{\partial^2}{\partial \x_i \partial \x_j}  \big( q_\syt(\x) \mathbf{D}_{i, j}(\syt, \x, q_\syt) \big),
    \quad q|_{\tau=0}(\x) = q_0(\x).
    \end{split}
\end{equation}
\Cref{eq:fokker-planck-vlasov2} is rewritten vectorially as 
\fi
\begin{equation}
\label{eq:fokker-planck-vlasov}
\begin{split}
 &   \frac{\partial q_\syt(\x)}{\partial \syt} = - \dive \Bigl( q_\syt(\x) \bigl[ \mathbf{F}(\syt, \x, q_\syt) - \mathbf{D}(\syt, \x, q_\syt) \nablax  \log q_\syt(\x) - \mathbf{d}(\syt, \x, q_\syt) \bigr] \Bigr), \\
 &   q|_{\tau=0}(\x) = q_0(\x),\quad \x \in \Omega \subseteq \Re^{\Nstate},
\end{split}
\end{equation}
where $\mathbf{D}(\syt, \x, q) \coloneqq (1/2)\sig(\syt, \x, q)\sig\tr(\syt, \x, q) \in \Re^{\Nstate \times \Nstate}$ is the diffusion tensor, 
$\dive = {\nabla}_{\x}\tr$ is the divergence operator, and $\mathbf{d}(\syt, \x, q) = \big(\dive \mathbf{D}(\syt, \x, q)\big)\tr \in \Re^{\Nstate}$. 

\begin{remark}
It is shown in \cite{Barbu_2010_FokkerPlanck} that \cref{eq:fokker-planck-vlasov} has solutions in the sense of distributions, and under general assumptions, $q_0 \in L^1(\Omega)$ implies $q_\tau \in C^0\big( [0,\infty), L^1(\Omega) \big)$. 
Existence of smooth solutions under more restrictive assumptions is discussed in \cite{Bouchut_1993_Vlasov-smooth,Grube_2023_Vlasov-strong,Degond_1986_Vlasov}. 
In \Cref{ass:connectivity} we make stronger smoothness assumptions on $q_\tau$.
\end{remark}

We seek to build particle flows \eqref{eq:ito-process-difussion} such that the stationary distribution of the corresponding Fokker-Planck-Vlasov equation \eqref{eq:fokker-planck-vlasov} is the posterior, i.e. $q_{\infty}(\x)=\PA(\x)$. 
Specifically, the KL divergence \cite{Kullback_1951_information} between $q_\syt(\x)$ and the target distribution $\PA(\x)$ is defined as 
\begin{equation}\label{eq:KL-div}
    \DKL (q_\syt\, \Vert \, \PA)  = \int_{\Omega} q_\syt(\x)\, \log \sfrac{q_\syt(\x)}{\PA(\x)}\, \dif \x, 
\end{equation}
where $\operatorname{support}(q_\syt) \subseteq \operatorname{support}(\PA) \coloneqq \Omega$. 
Our aim is to find drift $\mathbf{F}$ and diffusion $\sig$ terms such that the family of intermediate probability densities \eqref{eq:fokker-planck-vlasov}, initialized with $q_0(\x) \coloneqq \PB(\x)$, converges in KL-divergence to the posterior: 
\begin{equation}
\label{eq:KL-to-zero}
    \lim_{\syt\to\infty} \DKL \left(q_\syt\, \middle\Vert \, \PA \right) = 0.
\end{equation}
The process \eqref{eq:ito-process-difussion} provides an indexed family of intermediate random variables $\{\xsyt\}_{0 \leq \syt < \infty}$ that represents the dynamics of particles moving toward a sample of the target distribution.

We now present two theorems that help to define the KL divergence minimizing drift for \cref{eq:ito-process-difussion}.
Consider a smooth functional $\mathcal{A} : \Omega \times \Pi \to \Re^{\Nstate \times \Nstate}$ that maps each point in the domain and probability distribution to a symmetric positive definite  matrix $\mathcal{A}_\syt \coloneqq \mathcal{A}(\x, q_\syt)$, uniformly non-degenerate for any $q$.
Consider also the space of functions of finite second order scaled moments with respect to the probability density $q(\cdot)$ over $\Omega$, and define the following inner product:
\begin{subequations}
\begin{eqnarray}
\label{eq:q-dot-space}
L_{2,q,\mathcal{A}^{-1}}(\Omega) &\coloneqq& \Big\{ f : \Omega \to \Re^{\Nstate} \mid \int_\Omega q(\x)\,\norm{\mathcal{A}^{-1/2}(\x,q)\,f(\x)}^2 \, \dif\x < \infty\Big\}, \qquad \\
\label{eq:q-dot-product}
\langle f, g \rangle_{q,\mathcal{A}^{-1}} &\coloneqq&  \int_\Omega q(\x)\,f\tr(\x)\, \mathcal{A}^{-1}(\x,q)\,g(\x)\, \dif\x,
    \quad f,g \in L_{2,q,\mathcal{A}}(\Omega).
\end{eqnarray}
\end{subequations}
{
Thus, the drift is optimal in the sense that it is the direction of the largest decrease of the KL divergence between $q_\syt$ and $\PA$ in the inner product space defined by $\langle  \cdot, \cdot \rangle_{q_\syt, \mathcal{A}_\syt^{-1}}$. 
}

\begin{thm}\label{thm:optimal-drift}
Consider the process \cref{eq:ito-process-difussion}, and $(\syt,\x,q)  \mapsto \mathcal{A}_\syt \in \Re^{\Nstate \times \Nstate}$, a smooth functional that maps synthetic time, state, and probability densities to symmetric positive definite matrices.
The instantaneous optimal drift $\mathbf{F}$ that minimizes the KL-divergence \eqref{eq:KL-div} of the family of distributions governed by the Fokker-Planck equation \cref{eq:fokker-planck-vlasov} with respect to the dot-product $\langle \cdot, \cdot \rangle_{q_\syt,\mathcal{A}_\syt^{-1}}$ {in} \cref{eq:q-dot-product} is:
\begin{equation}\label{eq:optimal-drift}
\begin{split}
    \mathbf{F}(\syt, \x, q_\tau) 
    &= -\mathcal{A}_\syt(\x, q_\tau) \, \nablax  \log\sfrac{q_\syt(\x)}{\PA(\x)} +  \mathbf{D}(\syt, \x, q_\tau)\,\nablax \log{q_\syt(\x)} + \mathbf{d}(\syt, \x, q_\tau).
\end{split}
\end{equation}
The optimal drift \cref{eq:optimal-drift} depends on the current probability density $q_\syt$, as well as on diffusion term $\sig$ via $\mathbf{D}$ and $\mathbf{d}$.
    
The Fokker-Planck equation \cref{eq:fokker-planck-vlasov} under the optimal drift \cref{eq:optimal-drift} is:
\begin{equation}\label{eq:optimal-fokker-planck2}
    \frac{\partial q_\syt(\x)}{\partial \syt} =  \dive \left( q_\syt(\x)\, \mathcal{A}_\syt(\x, q_\tau) \, \nablax  \log \sfrac{q_\syt(\x)}{\PA(\x)} \right).
\end{equation}
\end{thm}

\begin{pf}
We omit all explicit arguments of the probability distributions and functions for brevity. The time derivative of the KL-divergence \cref{eq:KL-div} is:
\begin{equation*}
    \frac{\dif \DKL}{\dif \syt} = \int_{\Omega}\sfrac{\partial q_\syt}{\partial \syt} \left( \log \sfrac{q_\syt}{\PA} + 1\right) \dif \x,
\end{equation*}
and applying the Fokker-Planck equation \cref{eq:fokker-planck-vlasov} leads to:
\begin{equation*}
    \frac{\dif \DKL}{\dif \syt} = \int_{\Omega} - \dive  \left( q_\syt \bigl( \mathbf{F} - \mathbf{D} \nablax  \log q_\syt - \mathbf{d} \bigr) \middle)
    \middle( \log \sfrac{q_\syt}{\PA} + 1\right) \dif \x.
\end{equation*}
Since $\left. q_\syt \right|_{\partial \Omega} = 0$ by \Cref{ass:connectivity}, integrating by parts, and using \cref{eq:q-dot-product} leads to:
\begin{equation*}
    \frac{\dif \DKL}{\dif \syt} 
    = \left\langle  \mathbf{F} - \mathbf{D} \nablax  \log q_\syt - \mathbf{d}, \, \mathcal{A}_\syt \nablax  \left( \log \sfrac{q_\syt}{\PA}\right) \right\rangle_{q_\syt,\mathcal{A}_\syt^{-1}}.
\end{equation*}
Consequently the optimal drift $\mathbf{F}$ that maximizes the rate of decrease of the KL-divergence with respect to the dot product \cref{eq:q-dot-product} with scaling $\mathcal{A}_\syt$ is given by:
\begin{equation*}
    \mathbf{F}  - \mathbf{D} \nablax  \log q_\syt - \mathbf{d} = - \mathcal{A}_\syt \nablax  \left( \log \sfrac{q_\syt}{\PA}\right).
\end{equation*}
\end{pf}

\begin{thm}
The Fokker-Planck-Vlasov equation \cref{eq:optimal-fokker-planck2} evolves the probability distribution $q_\syt(\x)$ toward the unique steady-state $q_\infty(\x) = \PA(\x)$ a.e. (w.r.t. the Lebesgue measure), regardless of the initial condition $q_0(\x)$.
\end{thm}

\begin{pf}    
The proof follows \cite{Stuart_2020_gradient-EnKF}.
Consider the following modified Wasserstein distance between $\mu_0,\mu_1 \in \Pi$ \cite{Benamou_2000_CFD-transport}: 
\begin{equation}\label{eq:Wasserstein}
\begin{split}
    W^2(\mu_0,\mu_1) &= \Big\{\inf_{v_\syth} \int_0^1 \int_\Omega\, (v_\syth\tr\,\mathcal{A}_\syth^{-1}\,v_\syth)\,  q_\syth(\x)\,\dif \x \,\dif \syth 
    = \inf_{v_\syth} \int_0^1 \langle v_\syth, v_\syth \rangle_{q_\syth,\mathcal{A}_\syth^{-1}}\,\dif \syth \\
    & \qquad s.t. ~~\partial_\syth q_\syth = -\dive \left( q_\syth\,v_\syth \right), ~~ q_0 = \mu_0,~~q_1 = \mu_1 \Big\}.
\end{split}
\end{equation}
Define the following Riemannian metric tensor on the tangent space $T_{q_\syt} \Pi$ at $q_\syt$:
\begin{equation}\label{eq:Riemannian-metric}
    g_{q_\syt}( s_1, s_2 ) =   \int_\Omega\, (v_1\tr\,\mathcal{A}_\syt^{-1}\,v_2)\,  q_\syt(\x)\,\dif \x ~~ \textnormal{where} ~~ s_i = -  \dive \left( q_\syt\,v_i \right) \in T_{q_\syt} \Pi, ~ i=1,2,
\end{equation}
where $v_i = \nablax  \phi_{i}$, and $\phi_{i}$ are the unique solutions of the {linear} elliptic PDEs 
\begin{equation*}
- \dive ( q_\syt \nablax  \phi_{i} )= s_i, ~~ \x \in \Omega; \quad \phi_{i}\big|_{\x \in \partial \Omega} = 0,
\end{equation*}
and are vectors in the cotangent space, $\phi_{i} \in T^\ast_{q_\syt} \Pi$.
{(Note that $\int_\Omega s_i d\x = 0$.)}
The Riemannian gradient of the KL divergence $\mathcal{F}(q_\syt) \coloneqq \DKL \left(q_\syt\, \middle\Vert \, \PA \right)$, seen as a functional on $\Pi$, is defined as \cite{Stuart_2020_gradient-EnKF}:
\begin{equation}
\begin{split}
    g_{q_\syt}( \nabla_{q_\syt}\mathcal{F}, s ) &=   \int_\Omega\,\sfrac{\delta\mathcal{F}}{\delta q_\syt}(\x) \, s(\x)\,\dif \x = \int_{\Omega}\left( \log \sfrac{q_\syt}{\PA}(\x) + 1\right)\,s(\x)\,\dif \x \\
    &=   - \int_{\Omega}\left( \log \sfrac{q_\syt}{\PA}(\x) + 1\right)\,\dive \left( q_\syt\,v \right)\,\dif \x \\
    &=   \int_{\Omega} \left(\mathcal{A}_\syt\,\nablax  \log \sfrac{q_\syt}{\PA}(\x) \right)\tr\,\mathcal{A}_\syt^{-1}\,v(\x)\,  q_\syt (\x)\,\dif \x.
\end{split}
\end{equation}
 Comparing with the definition {in} \cref{eq:Riemannian-metric} of the Riemannian metric we obtain that
\begin{equation*}
    \nabla_{q_\syt}\DKL \left(q_\syt\, \middle\Vert \, \PA \right) = \dive \left( q_\syt\,\mathcal{A}_\syt\,\nablax  \log \sfrac{q_\syt}{\PA}(\x) \right),     
\end{equation*}
and that the FPE trajectory \cref{eq:optimal-fokker-planck2} under the optimal drift \cref{eq:optimal-drift} is a gradient flow for the KL divergence with respect to the modified Wasserstein distance \cref{eq:Wasserstein}:
\begin{equation}\label{eq:optimal-fokker-planck-gradient}
    \sfrac{\partial q_\syt(\x)}{\partial \syt} = -\nabla_{q_\syt}\DKL \left(q_\syt\, \middle\Vert \, \PA \right).
\end{equation}
Since $\DKL \left(q_\syt\, \middle\Vert \, \PA \right) \ge 0$ and it decreases monotonically along the FPE trajectory we conclude that \cref{eq:optimal-fokker-planck-gradient} converges to a steady state distribution $q_\infty(\x).$

A stationary distribution of the Fokker-Planck \cref{eq:optimal-fokker-planck2} satisfies:
\begin{equation}\label{eq:FP-stationary-distribution}
\begin{cases}
    \dive \left( q_\infty(\x)\, \mathcal{A}_\infty \, \nablax  \log\sfrac{q_\infty(\x)}{\PA(\x)}\right) = \dive \left( \PA(\x)\, \mathcal{A}_\infty \, \nablax  \sfrac{q_\infty(\x)}{\PA(\x)}\right) =0, & \x \in \Omega, \\
    q_\infty(\x) = 0, & \x \in \partial\Omega.
\end{cases}
\end{equation}
Since $\mathcal{A}_\infty>0$ by assumption, and $\PA(\x) > 0$ for $\x \in \Omega$, \cref{eq:FP-stationary-distribution} implies that $q_\infty(\x) \slash \PA(\x) =  const$ a.e. (w.r.t. the Lebesgue measure). 
Using the connectivity of $\Omega$, and since both distributions need to integrate to one, we conclude that $q_\infty(\x) = \PA(\x)$ a.e. 
\end{pf}

\begin{remark}
The optimal drift \eqref{eq:optimal-drift} consists of two terms, i.e., two forces acting on the particles. The term
\begin{equation}\label{eq:drift-part1}
    \mathcal{A}_\syt\, \big(\nablax  \log{\PA(\x)} - \nablax  \log{q_\syt(\x)}\big)
\end{equation}
is the scaled difference between the gradient-log-densities of the posterior and the intermediate distribution, and ensures that the  intermediate distribution is pushed toward the posterior one. 
The term
\begin{equation}\label{eq:drift-part2}
    \mathbf{D}(\syt, \x, q_\syt)\,\nablax  \log{q_\syt(\x)}   + \dive \mathbf{D}(\syt, \x, q_\syt),
\end{equation}
is an anti-diffusion term ``compensating'' for the stochastic term $\sig(\syt, \x)\dif \mathbf{W}_\syt$ of \cref{eq:ito-process-difussion}, and ensuring that the perturbations to any one realization of  intermediate variables still move towards being a realization of the analysis.
Deterministic dynamics are obtained by setting the diffusion to zero, $\sig = 0$. In this case the anti-diffusive force \cref{eq:drift-part2} is zero, and the optimal drift is given by the first term \cref{eq:drift-part1} only. 
The deterministic choice does not change the optimal FPE \cref{eq:optimal-fokker-planck2}.

\end{remark}

\begin{remark}\label{rem:general-func-der}
The KL-divergence \cref{eq:KL-div} is not the only way to quantify closeness to the posterior distribution.
Consider a general smooth finite functional on the space of smooth pdfs $\mathcal{F}: \Pi \to \Re$ and a smooth transform $g: \Re \to \Re$:
\begin{equation*}
\begin{split}
    \mathcal{F}(\pi) &\coloneqq g \left( \int_\Omega L(\x,\pi,\nablax  \pi)\, d\x \right), \quad \mathcal{F}(\pi) \ge 0 ~ \forall \pi, \\
    \sfrac{\delta \mathcal{F}(\pi)}{\delta \pi}(\x) &=  g' \left( \int_\Omega L(\x,\pi,\nablax  \pi)\, d\x \right)  \cdot \\
    &\quad \cdot \bigg( L_{\pi}\big(\x,\pi(\x),\nablax  \pi(\x)\big) - \operatorname{div} L_{\nablax  \pi}\big(\x,\pi(\x),\nablax  \pi(\x)\big) \bigg), \\
    \sfrac{\delta \mathcal{F}(\pi)}{\delta \pi}\Big|_{\pi = q}(\x) &= const~w.r.t.~\x
    \quad \Rightarrow \quad
    q(\x) = \PA(\x).
\end{split}
\end{equation*}
The optimal drift for decreasing $\mathcal{F}$ \cite{Reich_2021_FokkerPlanck} is:
\begin{equation*}
\begin{split}
    \mathbf{F}(\syt, \x, q_\syt) 
    &=  -\mathcal{A}_{\syt}\,\nablax \left(\sfrac{\delta \mathcal{F}(q_{\syt})}{\delta q_{\syt}}(\x)\right) + \mathbf{D}(\syt, \x, q_\syt)\, \nablax   \log q_{\syt}(\x) + \mathbf{d}(\syt, \x, q_\syt).
\end{split}
\end{equation*}
For example, consider the Renyi divergence functional with $\alpha > 0$, $\alpha \ne 1$:
\begin{equation*}
\begin{split}
    \mathcal{F}(q_{\syt}) &= D_{\alpha}(q_{\syt}\|\PA) =  \sfrac {1}{\alpha -1} \log \int_{\Omega} q_{\syt}(\x)^\alpha\,\PA(\x)^{1-\alpha}\,d\x, \\
    \mathbf{F}(\syt, \x, q_\tau) &= - c \, \mathcal{A}_\syt\,\nablax  \left( \sfrac{q_{\syt}(\x)}{\PA(\x)}\right)^{\alpha-1}  + \mathbf{D}(\syt, \x, q_\syt)\, \nablax   \log q_{\syt}(\x) + \mathbf{d}(\syt, \x, q_\syt),\\
    c &= \frac{\alpha}{\alpha - 1} \left( \int_{\Omega} q_{\syt}(\x)^\alpha\,\PA(\x)^{1-\alpha}\,d\x \right)^{-1}
\end{split}
\end{equation*}

\end{remark}

\subsection{Selection of the metric \texorpdfstring{$\mathcal{A}_\syt$}{Atau}}
\label{subsec:choosing-A}
%
The optimal drift {in} \cref{eq:optimal-drift} depends on the choice of $\mathcal{A}_\syt$, i.e., depends on the metric in which the minimum KL-divergence tendency is measured. Several special choices of $\mathcal{A}_\syt$ are discussed next.
\begin{enumerate}
\item The trivial choice, used in this paper, is $\mathcal{A}_\syt = \*I_{\Nstate}$, giving
\begin{equation}\label{eq:optimal-drift-trivial}
    \mathbf{F}(\syt, \x, q_\syt) = \nablax  \log{\PA(\x)} +  (\mathbf{D}(\syt, \x, q_\syt) - \*I_{\Nstate}) \nablax  \log{q_\syt(\x)} + \mathbf{d}(\syt, \x, q_\syt).
\end{equation}
The space \cref{eq:q-dot-space} are functions with finite second order moments.
\item In Stein variational gradient descent \cite{Liu_2016_Stein,Liu_2017_Stein} one chooses $\mathcal{A}_\syt = q_\syt(\x)\, \*I_{\Nstate}$ and $\sig = 0$, giving
\begin{equation}\label{eq:Stein-optimal-drift}
    \mathbf{F}(\syt, \x, q_\syt) = q_\syt(\x) \nablax  \log{\PA(\x)} - q_\syt(\x) \nablax  \log{q_\syt(\x)}.
\end{equation}
The space \cref{eq:q-dot-space} is $L_2(\Omega)$, the space of square integrable functions with respect to Lebesgue measure.
Equation \cref{eq:Stein-optimal-drift} can be seen as a scaling of the gradients obtained in \cref{eq:optimal-drift-trivial} without stochastic noise. 
Embedding the optimal drift \cref{eq:Stein-optimal-drift} in a RKHS with kernel $K(\cdot,\cdot)$ recovers the original formulation \cite{Liu_2016_Stein,Liu_2017_Stein}:
\begin{equation}
\label{eq:Stein-optimal-drift-RKHS}
\begin{split}
    \mathbf{F}(\syt, \x, q_\syt)\big|_{\textsc{rkhs}} &= \int_{\Omega} K(\x,\z)\,\mathbf{F}(\syt, \z)\, d\z  \\
    &\stackrel{ \cref{eq:Stein-optimal-drift} }{=} \expectw{\z \from q_\syt}{ K(\x,\z)\, \nablaz  \log{\PA(\z)} + \nablaz  K(\x,\z) }.
\end{split}
\end{equation}
\item The choice $\mathcal{A}_\syt = \mathbf{D}(\syt, \x, q_\tau)$ leads to {first-order, overdamped} Langevin dynamics \cite{Stuart_2020_gradient-EnKF,Reich_2019_interacting-Langevin}:
\begin{equation}
\label{eq:optimal-drift-Langevin}
    \mathbf{F}(\syt, \x, q_\syt) = \mathbf{D}(\syt, \x, q_\tau)\, \nablax  \log{\PA(\x)} + \mathbf{d} (\syt, \x, q_\tau),
\end{equation}
and this drift is optimal under the assumption that $\mathbf{D}(\syt, \x, q_\tau)$ has full rank. If $M = \Nstate$ and the drift matrix $\sig$ in \cref{eq:ito-process-difussion} is non-singular, then the space \cref{eq:q-dot-space} consists of functions $f$ for which $\sig^{-1}f$ has finite second order moments.
In the Langevin choice the optimal drift \cref{eq:optimal-drift-Langevin} does not depend on the current probability density $q_\syt(\x)$, and therefore can be computed very efficiently. The particles spread is ensured by diffusion only. 
\item An intuitive choice is the {sample covariance of the current distribution, i.e., $\mathcal{A}_\syt = \cov{q_\syt}$} \cite{Stuart_2020_gradient-EnKF}; in the context of Langevin dynamics \cref{eq:optimal-drift-Langevin} this choice leads to an affine-invariant flow \cite{Reich_2019_interacting-Langevin}.
\item Tapering $\mathcal{A}_\syt$ can be used to perform localization in a VFP context. The optimal drift acting on a certain state variable is restricted to depend only on gradient-log-density entries corresponding to ``nearby'' state variables. 
\end{enumerate}

\subsection{Discretization and parameterization of the optimal drift in the VFP filter}
\label{subsec:dicretization-and-parameterization}
\begin{table}[!ht]
\centering
\def\arraystretch{1.65}%
    \begin{tabular}{ |c|c| } 
    \hline
        Mixture & 
        \begin{tabular}{@{}c@{}}$\Prob(\x; \mathbf{\Theta}_{1:m}(\x) ) \propto \sum_{i=1}^{m} w_i\,\Prob_i(\x; \mathbf{\Theta}_i(\x))$,  \\[4pt]
        $\nablax  \log\Prob(\x; \mathbf{\Theta}_{1:m}(\x)) = \frac{\sum_{i=1}^{m}  w_i \, \big( \frac{\partial}{\partial \x}\Prob_i(\x; \mathbf{\Theta}_i(\x)) + \frac{\partial}{\partial \mathbf{\Theta}_i}\Prob_i(\x; \mathbf{\Theta}_i(\x))\frac{\partial \mathbf{\Theta}_i}{\partial \x} \big)}{\sum_{i=1}^{m} w_i \Prob_i(\x; \mathbf{\Theta}_i(\x))}$, \\[4pt]
        Simplifying assumption: $\partial \mathbf{\Theta}_i/\partial \x = 0,~ \forall i$.
        \end{tabular}\\
        \hline
        Kernel (K) & 
        \begin{tabular}{@{}c@{}}$\Prob(\x) \propto \frac{1}{\Nens}\sum_{i=1}^{\Nens} \!{K}(\x - \x_i)$,  \\
        $\nablax  \log\Prob(\x) =  \frac{\sum_{i=1}^{\Nens} \nablax  \!{K}(\x - \x_i)}{\sum_{i=1}^{\Nens} \!{K}(\x - \x_i)}$, \\
        $\!{K}$ is a positive definite kernel function.
        \end{tabular}\\
        \hline
        Gaussian (G) & 
        \begin{tabular}{@{}c@{}}$ \Prob(\x) \propto \exp(-\sfrac{1}{2}( \x - \xmean)\tr \P^{-1}\,( \x - \xmean))$,  \\
        $\nablax  \log\Prob(\x) = - \P^{-1}\,( \x - \xmean )$.
        \end{tabular}\\
        \hline
        Laplace (L) & \begin{tabular}{@{}c@{}}$\Prob(\x) \propto (\theta^\nu)\,\!{K}_\nu(\theta)$ ,\\
        $\nablax  \log\Prob(\x) = - \frac{2}{\theta}\frac{\!{K}_{\nu - 1}(\theta)}{\!{K}_\nu(\theta)}\P^{-1}\,( \x - \xmean )$, \\
        $\theta = \sqrt{2 ( \x - \xmean)\tr \P^{-1}\,( \x - \xmean)}, \quad \nu = 1-\Nstate/2$, \\ 
        $\!K_\nu$ is the modified Bessel function of the second kind~\cite{Abramovitz_1988_handbook,Kotz_2001_laplace}.
        \end{tabular}\\
        \hline
        Huber (H) & 
        \begin{tabular}{@{}c@{}}
        $   \nablax  \log\Prob(\x) = 
            \begin{cases}
                - \delta_1 \frac{2}{\theta}\frac{\!{K}_{\nu - 1}(\theta)}{\!{K}_\nu(\theta)}\P^{-1}\,(\x - \xmean) &\quad \delta_1 \frac{2}{\theta}\frac{\!{K}_{\nu-1}(\theta)}{\!{K}_\nu(\theta)} < \delta_2,\\
                -\delta_2\,\P^{-1}\, ( \x - \xmean ) &\quad \text{otherwise.}
            \end{cases}$        
        \end{tabular}\\
        \hline
        Cauchy (C) & 
        \begin{tabular}{@{}c@{}}
        \rule{0pt}{1ex} 
        $\Prob(\x) \propto \prod_{i=1}^n {\left[\pi\gamma_i\left(1 + \left(\frac{\x_i - \bar{\x}_i}{\gamma_i}\right)^2 \right) \right]}^{-1}$\\
        $\left[\nablax  \log\Prob(\x)\right]_i =-2\frac{\x_i-\bar{\x}_i}{\gamma^2_i + \left(\x_i-\bar{\x}_i\right)^2} $,\quad for $i=1,\dots,n$   
        \end{tabular}\\
        \hline
    \end{tabular}
    \caption{A collection of several parametrized probability distributions considered in this work, and the corresponding gradient-log-densities. The letters in the parentheses represent abbreviations of the distributions used to name the various families of VFP methods. For most of the distributions listed, the parameters are a semblance of centering $\xmean \in \Re^{\Nstate}$ (which may, but does not necessarily stand for the mean), and a semblance of spread $\*P \in \Re^{\Nstate \times \Nstate}$ (which may, but does not necessarily stand for covariance).}
    \label{tab:parameterized-families}
\end{table}
\def\arraystretch{1}

We now formulate the particle flow method discussed in \Cref{subsec:VFP-derivation} using a finite ensemble of particles.
Consider the ensemble $\Xenssyt \in \Re^{\Nstate \times \Nens}$ \eqref{eq:ensemble-definition} consisting of $\Nens$ particles $\x[e]_{\syt} \sim q_\syt(\cdot)$.
We refer to $\Xenssyt$ as the intermediate or current ensemble, and to $q_\syt$ as the intermediate or current distribution.
The McKean-Vlasov-It\^{o} process~\cref{eq:ito-process-difussion} acting on the ensemble is defined for each particle in artificial time as follows:
\begin{equation}\label{eq:ito-process-ensemble}
    \dif \xsyt^{[e]}  = \mathbf{F}(\syt, \xsyt^{[e]}, q_\tau)\, \dif \syt + \sig(\syt, \xsyt^{[e]}, q_\tau)\,\dif \*W_\syt, 
    \quad e = 1, \dots, \Nens.
\end{equation}
The optimal drift \cref{eq:optimal-drift-trivial} defined in \cref{thm:optimal-drift} acting on each particle \cref{eq:ito-process-ensemble} is:
\begin{equation}\label{eq:drift-ensemble}
    \mathbf{F}(\syt, \xsyt^{[e]}, q_\tau) =  {\mathcal{A}_\syt} \nablax \log{\PA(\xsyt^{[e]})} +  (\mathbf{D}(\syt, \xsyt^{[e]}) - {\mathcal{A}_\syt}) \nablax  \log{q_\syt(\xsyt^{[e]})} + \mathbf{d} (\syt, \xsyt^{[e]}),
\end{equation}
and depends both on the (continuous) analysis distribution $\PA$ and on the (continuous) intermediate distributions $q_\syt$, evaluated at the current particle state. 
Under the action of the flow, $\Xenssyt$ evolves toward an ensemble $\Xens_{\infty} = \Xensa$ of samples from the posterior distribution $\PA$ given by \cref{eq:Bayes-rule} or \cref{eq:Bayes-smooth}.

The drift term~\cref{eq:drift-ensemble}, requiring the gradient-log-likelihoods of the intermediate $q_\syt(\x)$ and the posterior $\PA(\x)$ probability densities, can be estimated in two ways. 
The first approach, proposed by Maoutsa et al.~\cite{Reich_2020_nabla-log-p}, expresses each analytical gradient-log-density $-\nablax  \log q_\syt(\x)$ and $-\nablax  \log\PA(\x)$ as the solution of a minimization problem. 
The second approach, employed in this paper, first reconstructs the continuous probability densities $q_\syt(\x)$ and $\PA(\x)$ using information from the ensembles $\Xenssyt$ and $\Xensb$, respectively, under appropriate assumptions. 
The corresponding gradient-log-densities are then evaluated on each particle. 
{Powerful kernelized dynamics can be obtained by embedding the drift $\mathbf{F}$ in an RKHS similar to~\cite{Pulido_2019_mapping-PF,Hu_2020_mapping-PF}, but without eliminating $q_\syt$.}

The VFP filter computes the optimal drift {in} \cref{eq:drift-ensemble} as follows:
\begin{enumerate}
    \item Assume the form of the prior distribution $\PB(\x)$, and fit the parameters of this distribution using the background ensemble $\Xensb$. Compute the corresponding negative gradient-log-likelihood function $-\nablax \log \PB(\x)$.
    \item By Bayes' rule~\cref{eq:Bayes-rule} the analysis gradient-log-likelihood is the sum of the gradient-log-likelihoods of the prior distribution and of the (known) observation error distribution:
    \begin{equation}\label{eq:analysis-gradient-loglikelihood}
    -\nablax \log{\PA}(\x) = -\nablax \log \PB(\x) - \nablax \log \PO(\x).
    \end{equation}
    \item Assume the form of the intermediate probability density $q_\syt(\x)$, and fit the parameters of this distribution using the current ensemble data $\Xenssyt$. Compute the corresponding negative gradient-log-likelihood function $-\nablax \log q_\syt(\x)$.
    \item Compute the optimal drift via formula \cref{eq:drift-ensemble} by evaluating the above gradients at particle states.
\end{enumerate}
The assumptions on the form of the prior distribution $\PB$ and of the intermediate distribution $q_\syt$ can differ from each other.
These choices dictate the parameterization of the VFP method.
We use the abbreviations in \cref{tab:parameterized-families} to distinguish between them.
For instance, we use the notation VFP(GH) to indicate a Gaussian assumption on the prior distribution and a Huber assumption on the intermediate distribution. 
VFPLn(G) is used to denote the Langevin variant of the Fokker-Planck dynamics with a Gaussian assumption on the prior.
\Cref{tab:parameterized-families} provides a non-exhaustive list of parameterized families of distributions: 
\begin{itemize}
    \item A general approach to parametrized distributions involves mixture modeling with an arbitrary set of parameters estimated from the corresponding ensemble. Though not  implemented in this work, it is of high research interest.
    \item The multivariate Gaussian distribution is an assumption similar to that made in ensemble Kalman filter methods. 
    \item The multivariate Laplace and the multivariate Huber distributions from {the field of} robust statistics~\cite{Huber_2011_Stats,Sandu_2017_robust-DA}.
    \item {In the VFP framework, the mapping particle filter~\cite{Pulido_2019_mapping-PF} and the high-dimensional flow filter~\cite{Hu_2020_mapping-PF} can be derived by embedding the optimal drift with $\mathcal{A}_\syt = q_\syt$ in an RKHS, eliminating $q_\syt$, and making kernel and Gaussian parameterizations on $\PB(\x)$.
    } 
\end{itemize}

\begin{remark}
The assumptions on $\PB(\x)$ and $\PO(\y \mid \x)$ directly impact the approximation of $\PA$, toward which particles converge in the limit. The assumptions and parameterizations made for $q_\syt$ change the way particles move towards the posterior, but not the limit of the process. 
\end{remark}

\begin{remark}
It is possible that the choices of parameterized families lead to intermediate distributions that do not converge in KL-divergence to the analysis~\cref{eq:KL-div}. 
For example, if $\PA$ is the product of a Gaussian ($\PO$) and Laplace ($\PB$) distributions, then neither a pure Gaussian nor a pure Laplace is a good assumption on $q_\syt$.
Thus, the KL divergence between $q_\syt$ and $\PA$ can then never be zero.
If we consider the parameterized family of intermediate distributions $\!Q$, then, instead of seeking a distribution such that the KL-divergence is zero~\cref{eq:KL-to-zero}, we instead aim to find the optimal analysis from the intermediate family that simply minimizes the KL-divergence:
\begin{equation*}
    \lim_{\syt\to\infty} \DKL \left(q_\syt\, \middle\Vert \, \PA\right) = \min_{q \in \!Q} \DKL \left(q\, \middle\Vert \, \PA\right).
\end{equation*}
\end{remark}

\subsection{Optimal drift in the VFP smoother}
\label{subsec:dicretization-and-parameterization-smoother}
Consider the {``strong-constraint''} smoothing posterior \cref{eq:Bayes-smooth} for a perfect model \cref{eq:model} with $\errm_i=0$ where $\x_{i} = \Model_{0,i}(\x_0)$ for any time $t_i$. 
Since model trajectories are fully determined by the initial conditions, the ``strong-constraint'' data assimilation \cite{Sandu_2011_assimilationOverview} is performed in the space of initial conditions. 
The ensemble of particles $\Xenssyt \in \Re^{\Nstate \times \Nens}$ \cref{eq:ensemble-definition} consists of $\Nens$ initial conditions $\x[e]_{0,\syt} \sim q_\syt(\x_0)$ that completely determine the $\Nens$ trajectories. (Here $0$ is the initial physical time, and $\tau$ is the synthetic time corresponding to changing particles/initial conditions.) 
The analysis probability density \cref{eq:Bayes-smooth} has the following gradient-log-likelihood:
\begin{equation}\label{eq:smoothing-analysis-gradient-loglikelihood}
\begin{split}
    \nabla_{\x_0}\log \PA(\x_0;\, \y_{0:K}) &= \nabla_{\x_0}\log \PB(\x_0) \\
    & \quad + \sum_{i=0}^K \Madj_{0,i}(\x_0)\,\nabla_{\x_i}\log \PO_{i}\big(\y_{i} \mid \x_i \big)\big|_{\x_i = \Model_{0,i}(\x_0)},
\end{split}
\end{equation}
where $\Madj_{0,i}(\x_0) \coloneqq \big( \dif \Model_{0,i}(\x_0) \slash \dif \x_0 \big)\tr$ is the adjoint model operator \cite{Sandu_2000_RosAdjoint,Sandu_2008_fdvarTexas}.

\begin{remark}[Traditional 4D-Var]
{Traditional 4D-Var  computes a maximum aposteriori estimate $\xa_{0}$ of the initial condition; the argument that maximizes the posterior probability density can be obtained by evolving the initial condition in synthetic time along the gradient-log-density \eqref{eq:smoothing-analysis-gradient-loglikelihood}:
\begin{equation}
\label{eqn:clasic-fdvar}
\begin{split}
   \sfrac{d}{d\tau} \x_{0,\tau} = \nabla_{\x_0}\log \PA(\x_{0,\tau};\, \y_{0:K}),
   \quad \x_{0,\tau} |_{\tau=0} = \xb_0, \quad
   \x_{0,\tau} \xrightarrow{\tau \to \infty} \xa_0.
\end{split}
\end{equation}
Evaluation of the posterior gradient-log-density values $\nabla_{\x_0}\log \PA$ requires one forward model run, followed by one adjoint model run.
}
\end{remark}

\begin{remark}[Ensemble of 4D-Vars]
\label{rem:En4dvar}
{The ``ensemble of 4D-Vars''  approach \cite{Lorenc_2012_nomenclature} performs $\Nens$ independent minimizations \eqref{eqn:clasic-fdvar} of negative log posterior densities corresponding to different samples of background states and perturbed observations $\y[e]$:
\begin{equation}\label{eq:en4dvar}
\begin{split}
&    \sfrac{d}{d\tau} \x[e]_{0,\tau}  = \nabla_{\x_0}\log \PA(\x[e]_{0,\tau};\, \y[e]_{0:K}), \\
    & \x[e]_{0,\tau} |_{\tau=0} = \xb[e]_0, \quad
   \x[e]_{0,\tau} \xrightarrow{\tau \to \infty} \xa[e]_0, \quad e = 1,\dots,\Nens.
\end{split}
\end{equation}
The $\Nens$ evolution equations \eqref{eq:en4dvar} are independent of each other (and possibly solved in parallel).
The result is an ensemble of analysis initial conditions $\xa[e]_0$, which samples the posterior exactly only when the posterior is Gaussian.}
\end{remark}

To apply the variational smoother, we compute the optimal drift  at each synthetic time $\syt$ according to the formula \cref{eq:drift-ensemble},
%
%
where, to simplify the discussion we consider here $\mathcal{A}_\syt = \Id$, and remove the explicit dependency on $\tau$. Each initial condition is then evolved in synthetic time using the stochastic differential equation \eqref{eq:ito-process-ensemble}:
\begin{equation}\label{eq:vfps}
\begin{split}
    \dif \x[e]_{0,\tau} &=  \nabla_{\x_0} \log{\PA(\x[e]_{0,\tau};\, \y_{0:K})}\, \dif \syt \\
    & \quad +  (\mathbf{D}(\x[e]_{0,\tau}) - \Id)\,  \nabla_{\x_0}  \log{q_\syt(\x[e]_{0,\tau})}\, \dif \syt + \mathbf{d} ( \x[e]_{0,\tau})\, \dif \syt \\
    &\quad + \sig(\syt, \xsyt^{[e]}, q_\tau)\,\dif \*W_\syt, \qquad
 e = 1, \dots, \Nens.
    \end{split}
\end{equation}

We use the name VFPS to refer to the {``strong-constraint''} VFP smoother \eqref{eq:vfps}. 
The computations of  gradient-log-density values are repeated for each synthetic time $\syt$, i.e., for each iteration of the underlying gradient-based minimization of the KL divergence \cref{eq:KL-to-zero}. 
Evaluation of the posterior gradient-log-density values $\nabla_{\x_0}\log \PA(\x[e]_0)$ requires $\Nens$ forward model runs, followed by $\Nens$ adjoint model runs, each started from a different initial condition $\x[e]_{0,\syt}$ (and all computed independently). A parametric approximation of $q_\syt(\x_0)$ is constructed as in the filtering case, and the corresponding gradient-log-density values $\nabla_{\x_0} \log{q_\syt(\x[e]_{0,\tau})}$ are calculated; this term uses all $\Nens$ initial conditions, but does not require additional model runs. 

\begin{remark}[VFPS]
\label{rem:vfps_vs_fdvar}
{
We compare the ``strong-constraint VFPS'' algorithm  \eqref{eq:vfps} with the ensemble of 4D-Vars approach \eqref{eq:en4dvar}.
\begin{itemize}
\item The ensemble of 4D-Vars approach \eqref{eq:en4dvar} runs $\Nens$ {\em independent} strong-constraint 4D-Var solutions \cite{Sandu_2011_assimilationOverview}. 
\item In contrast, for deterministic dynamics ($\sig=0$) the ``strong-constraint VFPS'' algorithm \eqref{eq:vfps} runs an ensemble of $\Nens$ {\em coupled} strong-constraint 4D-Var solutions:
\begin{equation*}
\begin{split}
    \sfrac{d}{d\tau} \x[e]_{0,\tau}  &=  \nabla_{\x_0} \log{\PA(\x[e]_{0,\tau};\, \y_{0:K})}
 - \nabla_{\x_0}  \log{q_\syt(\x[e]_{0,\tau})}, \quad
 e = 1, \dots, \Nens.
    \end{split}
\end{equation*}
The coupling is realized by the reconstructed current density $q_\syt(\x[e]_{0,\tau})$, which uses all $\Nens$ particles.
The VFPS solutions are stochastic for general $\sig \ne 0$.
\end{itemize} 
VFPS rigorously computes a sample from the posterior distribution. 
The coupling between particles makes VFPS different, and more rigorous (even for Gaussian posteriors ), than the ``ensemble of 4D-Vars'' sampling approach. A complete analysis of the differences between the two approaches is outside the scope of this work.}
\end{remark}

In case of an imperfect model \cref{eq:model} with model errors $\errm_i \from \Prob_i^{\Model}(\errm)$, {with the posterior \cref{eq:weak-Bayes-smooth}} the ``weak-constraint'' data assimilation \cite{Sandu_2011_assimilationOverview} is performed in the space of model trajectories. 
The ensemble of particles $\Xenssyt \in \Re^{\Nstate \times \Nens}$ \cref{eq:ensemble-definition} consists of $\Nens$ trajectories $\x[e]_{0:K,\syt} \sim q_\syt(\x_{0:K})$. 
Assuming that model errors at different times are independent of each other, and of observation errors, the analysis gradient-log-likelihood \cref{eq:Bayes-smooth} is:
\begin{equation}\label{eq:weak-smoothing-analysis-gradient-loglikelihood}
\begin{split}
    \nabla_{\x_i}\log \PA(\x_{0:K}) &= \nabla_{\x_i}\log \PB(\x_i) 
    + \nabla_{\x_i}\log \PO_{i}\big(\y_{i} \mid \x_i \big) \\
  &\quad + \nabla_{\errm}\log \Prob_{i-1}^{\Model}\big(\x_i - \Model_{i-1,i}(\x_{i-1})\big) \\
   &\quad- \Madj_{i,i+1}(\x_{i})\,\nabla_{\errm}\log \Prob_{i}^{\Model}\big(\x_{i+1} - \Model_{i,i+1}(\x_{i})\big).
\end{split}
\end{equation}
To compute the optimal drift \eqref{eq:drift-ensemble}, parametric approximations of the distribution $q_\syt(\x_{0:K})$ are needed. 
Under the approximation that current state probability densities at different physical times are independent, $q_\syt(\x_{0:K}) \approx \prod_{i=0}^{K} q_\syt(\x_i)$, one constructs parametric approximations of $q_\syt(\x_i)$ for each physical time, and computes $\nabla_{\x_i}\log q_\syt(\x[e]_i)$ as in the filtering case.
Evaluation of the posterior gradient-log-density values $\nabla_{\x_i}\log \PA(\x[e]_{0:K})$ for each physical time $i = 0, \dots, K$ requires computing solution differences $\x_{i+1} - \Model_{i,i+1}(\x_{i})$, and applying adjoint operators $\Madj_{i,i+1}(\x_{i})$. 
The state $\x[e]_{i,\syt}$ of each particle $e$ at each physical time $i$ moves under the flow of a different McKean-Vlasov-It\^{o} process \cref{eq:ito-process-ensemble}, which makes the entire computation highly parallelizable. 
{As in \cref{rem:En4dvar}, the ensemble of ``weak-constraint'' 4D-Var analyses with perturbed observations is exact only when the posterior is Gaussian.}

\subsection{Selection of the diffusion term}  
\label{subsec:diff}

To complete the description of the stochastic dynamics \cref{eq:ito-process-ensemble} one needs to select the diffusion term $\sig(\syt,\x)$. 
Recall that the optimal drift \cref{eq:optimal-drift} depends on the diffusion term, however the resulting Fokker-Planck equation \cref{eq:optimal-fokker-planck2} does not. 
Nevertheless, the choice of $\sig(\syt,\x)$ does impact the implementation of the algorithm, as well as its practical performance, given the finite number of particles and the different approximations made when reconstructing probability densities.
The trivial choice $\sig(\syt,\x) \equiv 0$ can be made to ensure deterministic particle dynamics. Using the optimal drift \cref{eq:optimal-drift} the process \cref{eq:ito-process-ensemble} becomes:
\begin{equation}\label{eq:deterministic-process-ensemble}
    \sfrac{\dif}{\dif \syt} \xsyt^{[e]}  = - \mathcal{A}_\syt\, \nablax  \log\sfrac{q_\syt(\xsyt^{[e]})}{\PA(\xsyt^{[e]})},\quad \x[e]_0 = \xb[e],
    \quad e = 1, \dots, \Nens.
\end{equation}
Our experience with running the experiments indicate that, during the first data assimilation cycles, the deterministic VFP method \eqref{eq:deterministic-process-ensemble} successfully transforms particles from background samples into analysis samples. 
However, after multiple assimilation cycles, the performance of the filter deteriorates considerably due to the phenomenon of particle collapse.
{While deterministic dynamics are exact for an infinite number of particles, they lead to biased analysis for finite sizes, and ultimately lead to particle collapse.}
Using stochastic dynamics, i.e., a non-zero diffusion $\sig(\syt,\x)$, is akin to performing rejuvenation in particle filters \cite{popovamit, Reich_2013_ETPF, Reich_2017_ETPF_SOA}, and presents a natural approach to alleviate this problem.
{In our view, stochastic dynamics alleviates, to a large extent, the bias issue inherent with a finite number of particles.}

Since particles are physical model states, the choice of the diffusion term should ensure that the stochastic perturbations $\sig(\syt, \x)\dif \mathbf{W}_\syt$ do not push the particle states outside physical regimes.
These perturbations should respect the scaling of different components, correlations between variables, and the quasi-equilibria of the system.
To this end, a reasonable choice is to use a scaled square root of  the forecast particle covariance ($\sig(\syt, \x) = \diffu\, (\Pb)^{1/2}$), or of the current particle covariance  ($\sig(\syt, \x) = \diffu\, \cov{\Xenssyt}^{1/2}$) \cite{Stuart_2020_gradient-EnKF}, or of a climatological covariance ($\sig(\syt, \x) = \diffu\, {\*B}^{1/2}$), {where the climatological covariance $\*B$ is a data driven estimate of the covariance, typically the autocovariance, and $\diffu$ is a scaling parameter}.
An approximation to the square root of the forecast covariance is given by the scaled ensemble anomalies ($(\Pb)^{1/2} = (\Nens-1)^{-1/2} (\Xensb - \xmeanb)$), in which case the stochastic perturbations are similar to the rejuvenation done in ETPF \cite{Reich_2019_interacting-Langevin}.
A similar approximation can be employed for the square root of the current particle covariance. Note that choosing the diffusion matrix $\sig(\syt, \x, q)$ to be independent of the states, makes $\mathbf{d} = \dive \mathbf{D} = 0$, and simplifies the drift computation.

\subsection{Regularization of the particle flow}
\label{subsec:reg}
With a finite number of particles in a large state space, particle and ensemble filters may suffer from ensemble collapse.
As discussed in \Cref{subsec:diff}, the stochastic diffusion in VFP plays a role similar to that of particle rejuvenation in traditional particle filters.
However, since the diffusion is dependent on many parameters, there is no guaranteed prevention of ensemble collapse.
For this, we consider a regularization of the particle flow, which adds an additional drift term to the dynamics in \cref{eq:ito-process-difussion}, that favors particle spread in the state space. 
Thus, the resulting drift has one component that pushes the particles toward a sample of the posterior, by minimizing the KL divergence {in} \cref{eq:KL-div}, and a regularization component that pushes the particles apart.

\subsubsection{General considerations and particle interaction}
{We now ask the question: what is the effect of a finite number of particles on the optimal dynamics \eqref{eq:optimal-drift}?
Assuming that $q_\tau$ is known exactly, the optimal drift \eqref{eq:optimal-drift} applied to each particle (drifts differ due to different particle states $\x_\tau$) ensures that the marginal distribution of each particle approaches the posterior, $q_\infty = \mathcal{P}^a$.
But what is the joint distribution of all particles? 
Are converged particle states $\{\x[e]_\infty\}_{1 \le e \le \Nens}$ {\it independent} samples from the posterior distribution, as desired? 
When $q_\tau$ is not known but is approximated from all available particle states $\{\x[e]_\tau\}_{1 \le e \le \Nens}$, the dynamics of each particle depends on all other particles, and therefore the particle states are correlated random variables.
}	
	 
{
In order to answer this question, motivated by \cite{Reich_2021_FokkerPlanck}, we consider an ensemble of particles $\Xens \in \Re^{\Nstate \times \Nens}$ \cref{eq:ensemble-definition} as one large system of interacting particles.
Formally, we define one large state vector that stacks all particles as follows:
\begin{equation}\label{eq:reich-coupling}
    \mathfrak{X} \coloneqq \left[ \x^{[1]\mathrm{T}}, \x^{[2]\mathrm{T}}, \cdots, \x^{[\Nens]\mathrm{T}} \right]\tr \in \Re^{\Nens\Nstate \times 1}, \quad
    \mathfrak{X}^{[e]} \coloneqq \x[e]
\end{equation}
where $1 \le e \le \Nens$.
The current and the background vectors \eqref{eq:reich-coupling} of interacting particles  are denoted by $\mathfrak{X}_\syt$ and $\mathfrak{X}^{\rm b}$, respectively.}
{
\begin{remark}
The methodology discussed in \cref{subsec:dicretization-and-parameterization} and \cref{subsec:dicretization-and-parameterization-smoother} uses the finite number of particles to build build parametric approximations of $\PA$ and $q_\syt$.
With some abuse of notation, let $\PA(\x,\mathfrak{X}^{\rm b})$, and $q_\syt(\x,\mathfrak{X}_\syt)$ also denote these reconstructed probabilities, where their dependency on the ensembles used to fit parameters is made explicit. 
%
%
The general drift $\mathbf{F}(\syt, \xsyt^{[e]}, q_\syt)$ and diffusion $\boldsymbol{\sigma}(\syt, \xsyt^{[e]}, q_\syt)$ terms that act on particle $e$ are represented in parametric form as $\mathbf{F}(\syt, \xsyt^{[e]},\mathfrak{X}_\syt)$ and $\boldsymbol{\sigma}(\syt, \xsyt^{[e]},\mathfrak{X}_\syt)$, respectively. 
Consequently, the evolution of $e$-th particle depends on the states of all particles, which results in a coupling of particle dynamics.
\end{remark}
}

{
To perform variational filtering using $\Nens$ particles, we evolve $\mathfrak{X}_\syt$ through the following stochastic dynamics:
\begin{equation}\label{eq:ito-step-system}
    \dif \mathfrak{X}_{\syt} = \mathfrak{F}(\syt, \mathfrak{X}_{\syt})\, \dif \syt + \mathfrak{S}(\syt, \mathfrak{X}_{\syt})\, \dif \mathfrak{W}_\syt,
\end{equation}
where $\mathfrak{F}: \Re_+ \times \Re^{\Nstate\Nens} \to \Re^{\Nstate\Nens}$ is the optimal drift, $\mathfrak{S}: \Re_+ \times  \Re^{\Nstate\Nens} \to \Re^{\Nstate\Nens \times M\Nens}$ is the diffusion, and $\mathfrak{W} \syt \in \Re^{M\Nens}$ is a Wiener process. 
A sequence of random variables whose joint distribution is invariant under any reordering is called exchangeable ~\cite{ONeill_2009_exchangeability}.
This concept is weaker than independent and identically distributed (i.i.d) as all i.i.d sequences are trivially exchangeable, while the vice versa is not true. 
Let $\mathfrak{Q}_\syt(\mathfrak{X})$ be the joint probability density of $\mathfrak{X}_\syt$ on $\Re^{\Nens\Nstate}$ evolving under the process \cref{eq:ito-step-system}. 
Assuming particle exchangeability, the joint probability distribution $\mathfrak{Q}$ is independent of the particle stacking order in \cref{eq:reich-coupling}.
Therefore, as the marginal probability densities of all particles are equal to each other ($q_\syt$), we can define the joint density as 
\begin{equation}\label{eq:coupled-copula}
    \mathfrak{Q}_\syt(\mathfrak{X}) \coloneqq r_\syt(\mathfrak{X})  \prod_{e=1}^{\Nens} q_\syt(\x[e], \mathfrak{X}_\syt),
\end{equation}
where $r_\syt$ couples the marginal densities of individual particles. 
As the goal is to push particles toward i.i.d. samples from the posterior, the target probability density on $\Re^{\Nens\Nstate}$ is $\PA(\mathfrak{X}) = \prod_{e=1}^{\Nens} \PA(\x[e], \mathfrak{X}^{\rm b})$. 
The optimal drift \cref{eq:optimal-drift} for the system given in \cref{eq:ito-step-system} is 
\begin{equation}\label{eq:optimal-drift-full-system}
    \begin{split}
        \mathfrak{F}(\syt, \mathfrak{X}, \mathfrak{Q}_\syt) 
        &= -\mathfrak{A}_\syt(\mathfrak{X}, \mathfrak{Q}_\syt) \, \nabla_{\mathfrak{X}} \log \sfrac{\mathfrak{Q}_\syt(\mathfrak{X})}{\PA(\mathfrak{X})}\\ 
        &+ \mathfrak{D}(\syt, \mathfrak{X}, \mathfrak{Q}_\syt)\,\nabla_{\mathfrak{X}} \log{\mathfrak{Q}_\syt(\mathfrak{X})} + \mathfrak{d}(\syt, \mathfrak{X}, \mathfrak{Q}_\syt),
    \end{split}
\end{equation}
where $\mathfrak{D} = \frac{\mathfrak{S} \mathfrak{S}^\mathrm{T}}{2}$, and $\mathfrak{d} = \nabla_{\mathfrak{X}}^\mathrm{T} \mathfrak{D}$ (similar to definitions of $\mathbf{D}$ and $\mathbf{d}$ in \cref{subsec:VFP-derivation}).
Consider the optimal drift for the system of particles in \cref{eq:optimal-drift-full-system} along with simplifying assumptions such as $\mathfrak{S}(\syt, \mathfrak{X}_\syt) \coloneqq \operatorname{blkdiag}_{e = 1 \dots \Nens} \left\{ \sig(\syt, \x[e]_\syt,  q_\syt) \right\}$, and $\mathfrak{A}_\syt(\mathfrak{X}, \mathfrak{Q}_\syt) \coloneqq \operatorname{blkdiag}_{e = 1 \dots \Nens} \left \{ \mathcal{A}_\syt(\x,  q_\syt) \right\}$; the drift component acting on each particle $e$ is:
\begin{equation}\label{eq:optimal-drift-particle-system}
\begin{split}
    \widehat{\mathbf{F}}(\syt, \x[e]_\syt,\mathfrak{X}_\syt) &= \mathcal{A}_\syt \nabla_{\x[e]_\syt} \log{\PA(\x[e]_\syt,\mathfrak{X}^{\rm b})} + \mathbf{d}(\syt, \x^{[e]}_\syt, \mathfrak{X}_\syt) \\
    &+ (\mathbf{D}(\syt, \x[e]_\syt,\mathfrak{X}_\syt) - \mathcal{A}_\syt) \, \left( \nabla_{\x[e]_\syt} \log{ \left(q_\syt(\x[e]_\syt,\mathfrak{X}_\syt) r_\syt(\mathfrak{X}_\syt) \right)}  \right),\\
    &= \mathbf{F}(\syt, \x[e]_\syt,\mathfrak{X}_\syt)
    + \left( \mathbf{D}(\syt, \x[e]_\syt,\mathfrak{X}_\syt) - \mathcal{A}_\syt \right) \, \left( \nabla_{\mathfrak{X}^{[e]}_\syt} \log r_\syt(\mathfrak{X}_\syt) \right)\\ 
    &+ \left( \mathbf{D}(\syt, \x[e]_\syt,\mathfrak{X}_\syt) - \mathcal{A}_\syt \right) \left( \nabla_{\mathfrak{X}^{[e]}_\syt} \log q_\syt(\x[e]_\syt,\mathfrak{X}_\syt) \right).
\end{split}
\end{equation}
%
The optimal drift for the interactive particle system \cref{eq:optimal-drift-particle-system} consists of the optimal drift for an individual particle \cref{eq:drift-ensemble}, and an additional term --- the gradient-log of the probability coupling term \cref{eq:coupled-copula} $\nabla_{\x[e]_\syt} \log r_\syt(\mathfrak{X}_\syt)$ and the gradient-log of intermediate density with respect to the parametrization $\nabla_{\mathfrak{X}^{[e]}_\syt} \log q_\syt(\x[e]_\syt,\mathfrak{X}_\syt)$ --- that ``correct'' for the finite number of particles.
Thus, we show the existence of gradient-log coupling term that is missing in the ensemble of particle of framework as in \cref{eq:drift-ensemble}. 
\begin{remark}
    For these simplified assumptions, Langevin dynamics will remain unbiased as $\*D = \mathcal{A}_\syt$.
\end{remark}
For a more general case, without any simplifying assumptions on $\mathfrak{S}(\syt, \mathfrak{X}_\syt)$, and $\mathfrak{A}_\syt(\mathfrak{X}, \mathfrak{Q}_\syt)$, we have additional terms as
\begin{equation}\label{eq:optimal-drift-particle-system-extra}
\begin{split}
    \widehat{\mathbf{F}}(\syt, \x[e]_\syt,\mathfrak{X}_\syt) &= \mathbf{F}(\syt, \x[e]_\syt,\mathfrak{X}_\syt) + \sum_{\substack{i = 1 \\ i \neq e}}^{\Nens} \mathcal{A}^{[e,i]}_\syt \nabla_{\x[i]_\syt} \log{\PA(\x[i]_\syt,\mathfrak{X}^{\rm b})}\\ 
    &+ \sum_{\substack{i = 1 \\ i \neq e}}^{\Nens} \operatorname{div}_{\x[i]_\syt} \left( \mathbf{D}^{[e,i]}(\syt, \mathfrak{X}_\syt) \right) \\
    &+ \sum_{\substack{i = 1 \\ i \neq e}}^{\Nens} \left( \mathbf{D}^{[e,i]}(\syt, \mathfrak{X}_\syt) - \mathcal{A}^{[e,i]}_\syt \right) \, \left( \nabla_{\x[i]_\syt} \log q_\syt(\x[i]_\syt,\mathfrak{X}_\syt)  \right)\\ 
    &+ \sum_{i = 1}^{\Nens} \left( \mathbf{D}^{[e,i]}(\syt,\mathfrak{X}_\syt) - \mathcal{A}^{[e,i]}_\syt \right) \, \left( \nabla_{\mathfrak{X}^{[i]}_\syt} \log r_\syt(\mathfrak{X}_\syt) \right)\\ 
    &+ \sum_{i = 1}^{\Nens} \left( \mathbf{D}^{[e,i]}(\syt,\mathfrak{X}_\syt) - \mathcal{A}^{[e,i]}_\syt \right) \left( \nabla_{\mathfrak{X}^{[i]}_\syt} \log q_\syt(\x[i]_\syt,\mathfrak{X}_\syt) \right),
\end{split}
\end{equation}
where $\mathbf{D}^{[e,i]}(\syt,\mathfrak{X}_\syt)$ is the $(e, i)$-th block (in $\mathbb{R}^{n \times n}$) of $\mathfrak{D}(\syt,\mathfrak{X}_\syt)$, and $\mathcal{A}^{[e,i]}_\syt$ is the $(e, i)$-th block (in $\mathbb{R}^{n \times n}$) of $\mathfrak{A}_\syt(\mathfrak{X}, \mathfrak{Q}_\syt)$.
%
}

\begin{example}
Consider the case of first-order, overdamped Langevin dynamics where $\mathcal{A}_\syt = \mathbf{D}(\syt, \mathfrak{X}_\syt)$ is the empirical covariance of the current ensemble.
We have (see also \cite{Reich_2019_interacting-Langevin}):
\begin{equation*}
\begin{split}
    \mathbf{D}(\syt, \mathfrak{X}_\syt) &= \sfrac{1}{\Nens-1} \sum_{i=1}^{\Nens} \big(\x[i]_\syt - {\textstyle \sum_{j=1}^{\Nens}} \x[j]_\syt/\Nens \big)\,\big(\x[i]_\syt - {\textstyle \sum_{j=1}^{\Nens}} \x[j]_\syt/\Nens \big)\tr, \\
    \nabla_{\x[e]_\syt}\tr \mathbf{D}(\syt, \mathfrak{X}_\syt) 
    &= \sfrac{1}{2(\Nens-1)}  \big(\x[e]_\syt - {\textstyle \sum_{j=1}^{\Nens} \x[j]_\syt /\Nens}  \big), \\
    \widehat{\mathbf{F}}(\syt, \x[e]_\syt,\mathfrak{X}_\syt) & = \mathbf{F}(\syt, \x[e]_\syt,\mathfrak{X}_\syt) +  \sfrac{1}{2(\Nens-1)}  \big(\x[e]_\syt - {\textstyle \sum_{j=1}^{\Nens} \x[j]_\syt /\Nens}  \big).
\end{split}
\end{equation*}
\end{example}

\begin{example}
Consider the case of a multivariate Gaussian probability density with known statistics (no empirical statistics are used for parameterization, and therefore the corresponding derivatives with respect to parametrization ensembles are zero):
\begin{equation}\label{eq:gaussian-full-state}
    \mathfrak{Q}_\syt(\mathfrak{X}) 
    = \mathcal{N}\left(\mathfrak{X} \mid \bar{\mathfrak{X}}_\syt,  
    \begin{bmatrix}
        \P_\syt & \cdots & \widehat{\P}_\syt \\
        \vdots & \ddots & \vdots \\
        \widehat{\P}_\syt & \cdots & \P_\syt 
    \end{bmatrix}
    \right), \quad \bar{\mathfrak{X}}_\syt = 
    \begin{bmatrix}
        \xmean_\syt\\
        \vdots  \\
        \xmean_\syt 
    \end{bmatrix}.
\end{equation}
All particle means $\xmean_\syt$, covariances $\P_\syt$, and cross-covariances $\widehat{\P}_\syt$ are equal due to exchangeability. Marginalizing \cref{eq:gaussian-full-state} leads to the density of each particle: 
\begin{equation*}
q_\syt(\x) = \mathcal{N}\big(  \x \mid \xmean_\syt, \P_\syt \big).
\end{equation*}
The coupling term \cref{eq:coupled-copula} reads:
\begin{equation}\label{eq:gaussian-copula-forces}
    r_\syt(\mathfrak{X}) \propto \exp \left( -\sfrac{1}{2} \left(\mathfrak{X} - \bar{\mathfrak{X}}_\syt \right)\tr \mathfrak{C}_\syt (\mathfrak{X} - \bar{\mathfrak{X}}_\syt) \right)
\end{equation}
where the precision matrix $\mathfrak{C}_\syt$ is given by 
\begin{equation}
    \mathfrak{C}_\syt = 
    \begin{bmatrix}
        \C_\syt & \cdots & \widehat{\C}_\syt \\
        \vdots & \ddots & \vdots \\
        \widehat{\C}_\syt & \cdots & \C_\syt \\
    \end{bmatrix}
    =
    \left(
    \begin{bmatrix}
        \P_\syt & \cdots & \widehat{\P}_\syt \\
        \vdots & \ddots & \vdots \\
        \widehat{\P}_\syt & \cdots & \P_\syt \\
    \end{bmatrix}^{-1} - 
    \begin{bmatrix}
        \P_\syt & \cdots & \*0 \\
        \vdots & \ddots & \vdots \\
        \*0 & \cdots & \P_\syt \\
    \end{bmatrix}^{-1}
    \right),
\end{equation}
with elements 
\begin{equation*}
\begin{split}
    \widetilde{\C}_\syt &\coloneqq\P_\syt^{-1}\, \widehat{\P}_\syt\,\left(\sfrac{1}{\Nens - 1}\,\P_\syt  + \sfrac{\Nens - 2}{\Nens-1} \widehat{\P}_\syt - \widehat{\P}_\syt \P_\syt^{-1} \widehat{\P}_\syt \right)^{-1}, \\
    \C_\syt 
    & = \P_\syt^{-1}\, \widehat{\P}_\syt\,\widetilde{\C}_\syt, \qquad
     \widehat{\C}_\syt 
    =  -\sfrac{1}{\Nens - 1}\,\widetilde{\C}_\syt.
\end{split}
\end{equation*}
If the particles are independent, i.e., $\widehat{\P}_\syt = 0$ in \cref{eq:gaussian-full-state}, then $\log r_\syt \equiv 0$ \cref{eq:gaussian-copula-forces}.
\ifreport
Simplifying \cref{eq:gaussian-copula-forces} leads to
\begin{equation*}
    \log r_\syt(\mathfrak{X}) = -\sfrac{1}{2} \sum_{e = 1}^{\Nens} \left( (\x[e]_\syt - \xmean_\syt)\tr \C_\syt (\x[e]_\syt - \xmean_\syt) + \sum_{\substack{i = 1 \\ i \neq e}}^{\Nens} (\x[e]_\syt - \xmean_\syt)\tr \widehat{\C}_\syt (\x[i]_\syt - \xmean_\syt) \right) + c.
\end{equation*}
\fi
The log gradient of the coupling term \cref{eq:gaussian-copula-forces} for each particle is given by 
\begin{equation*}
\begin{split}
    &\nabla_{\x[e]_\syt} \log r_\syt(\mathfrak{X}) 
    = - \C_\syt (\x[e]_\syt - \xmean_\syt) -\sfrac{1}{2} \widehat{\C}_\syt  \sum_{\substack{i = 1 \\ i \neq e}}^{\Nens} (\x[i]_\syt - \xmean_\syt) \\
    &\qquad = - \left( \C_\syt -\sfrac{1}{2} \widehat{\C}_\syt \right) (\x[e]_\syt - \xmean_\syt) + \sfrac{\Nens}{2(\Nens - 1)}\,\widetilde{\C}_\syt \,\left( {\textstyle \sum_{i = 1}^{\Nens}  \x[i]_\syt /\Nens} - \xmean_\syt \right).
\end{split}
\end{equation*}
The first term is a force that pushes the particle away from the ensemble mean, therefore favoring ensemble spread. 
The second term, applied equally to all particles, is the random sampling error for particle mean, scaled by a factor that remains bounded for $\Nens > 1$ (note that $\widetilde{\C}_\syt$ is bounded for $\Nens > 1$). 
\end{example}

The {term $\nabla_{\x[e]_\syt} \log r_\syt(\mathfrak{X}_\syt)$ in} \cref{eq:optimal-drift-particle-system} that corrects for a finite number of particles while maintaining the ensemble spread is difficult to estimate.
For this reason, we consider modeling the interaction between particles via an interaction potential, and modifying the optimal drift such as to ensure particle independence and hence, maintain its spread.

\subsubsection{Modeling the particle interaction and regularization}

Let $\kappa : \Re^{\Nstate \times \Nstate} \to \Re$ be a smooth potential function that models the interaction between particles (specifically, $\kappa(\x,\widehat{\x})$ represents the interaction potential between $\x$ and $\widehat{\x}$ {that are assumed non-independent}).
We add to the  KL divergence functional a regularization term given by the average potential:
\begin{equation}\label{eq:regularized-KL-divergence}
    \DhKL(q_\syt \,\Vert\, \PA) = \DKL(q_\syt \,\Vert\, \PA) + \beta_\syt\, I_\syt, \quad I_\syt = \expectw{\x,\widehat{\x} \from q_\syt}{\kappa(\x,\widehat{\x})},
\end{equation} 
where the parameter $\beta_\syt$ is a pseudo-time dependent scalar that determines the strength of the regularization term. 
Minimization of \cref{eq:regularized-KL-divergence} decreases $\DKL(q_\syt \,\Vert\, \PA)$, therefore pushes particles toward the posterior, but also decreases the interaction between particles.

\begin{example}[Mutual information]
Mutual information, a non-negative real number, seeks to measure the independence of two random variables with zero indicating independence and any other value indicating {the} degree of dependence. 
We assume that the coupling \cref{eq:coupled-copula} of $\Nens$ exchangeable random variables is well described by a \textit{smooth} coupling between pairs of random variables $r_\syt(\x[j],\x[e])$ with $j \ne e$. 
Let two distinct particles with joint probability $\mathfrak{Q}_\syt(\cdot, \cdot)$, marginals $q_\syt(\cdot)$, and coupling term $r_\syt(\cdot,\cdot)$ \cref{eq:coupled-copula}.
The mutual information \cite{Cover_1999_infotheorybook} between the two distinct particles (random variables) is given by 
\begin{equation*}
\begin{split}
    I_\syt & = -\int_\Omega \int_\Omega \mathfrak{Q}_\syt(\x, \hat{\x}) \log{\left( \frac{\mathfrak{Q}_\syt(\x, \hat{\x})}{q_\syt(\x)q_\syt(\hat{\x})} \right)} \,\dif \x \, \dif \hat{\x} 
    =  \expectw{\x,\widehat{\x} \from q_\syt}{\kappa(\x,\widehat{\x})}, \\
     \kappa(\x,\widehat{\x}) &=  - r_\syt(\x,\widehat{\x})\, \log r_\syt(\x,\widehat{\x}).
\end{split}
\end{equation*}
Minimizing the regularized KL divergence \cref{eq:regularized-KL-divergence} pushes particles toward the posterior, but also nudges particles toward independence.
\end{example}

Following the same reasoning as in \cref{thm:optimal-drift}, the optimal drift $\widehat{\mathbf{F}}$ that minimizes the regularized KL divergence \cref{eq:regularized-KL-divergence} is given by 
\begin{equation}
\label{eq:optimal-drift-regularized}
    \widehat{\mathbf{F}}(\syt, \x, q_\syt) = \mathbf{F}(\syt, \x, q_\syt) - \beta_\syt\, \mathcal{A}_\syt\, \expectw{\hat{\x} \from q_\syt} { \nablax \kappa_\syt(\x, \hat{\x}) },
\end{equation} 
where $\mathbf{F}(\syt, \x, q_\syt)$ is the optimal drift without regularization \cref{eq:optimal-drift}.

\begin{remark}
Discretization of \cref{eq:optimal-drift-regularized} using $\Nens$ particles leads to the following drift acting on each particle $e = 1, \dots, \Nens$: 
\begin{equation}
\label{eq:regularized-drift}
    \widehat{\mathbf{F}}(\syt, \xsyt^{[e]},\mathfrak{X}_\syt) = \mathbf{F}(\syt, \xsyt^{[e]},\mathfrak{X}_\syt) - \sfrac{\beta_\syt}{\Nens} \,\mathcal{A}_\syt\, \sum_{\substack{i = 1 \\ i \neq e}}^{\Nens}  \nabla_{\xsyt^{[e]}} \kappa_\syt(\xsyt^{[e]}, \xsyt^{[i]}).
\end{equation}
The negative gradient $-\nabla_{\xsyt^{[e]}} \kappa(\xsyt^{[e]}, \xsyt^{[i]})$ in \cref{eq:regularized-drift} can be viewed as the corresponding repelling regularization force between particles $e$ and $i$, exerted on particle $e$; we do not include a repelling force of a particle on itself. 
This also holds true in the Gaussian case, as seen in \cref{eq:gaussian-copula-forces}. 
\end{remark}

The Fokker-Planck-Vlasov equation under the regularized drift \cref{eq:optimal-drift-regularized} is:
\begin{equation}\label{eq:optimal-fokker-planck-regularized}
    \frac{\partial q_\syt(\x)}{\partial \syt} = - \dive \left( q_\syt(\x)\, \mathcal{A}_\syt \, \nablax  \big(  \log{\sfrac{\PA(\x)}{q_\syt(\x)}}
    - \beta_\syt\, \expectw{\hat{\x} \from q_\syt}{\kappa_\syt(\x, \hat{\x})} \big) \right),
\end{equation}
and its stationary distributions are characterized by
\begin{equation*}
    \log{\sfrac{\PA(\x)}{q_\infty(\x)}}
    - \beta_\infty\, \expectw{\hat{\x} \from q_\infty}{\kappa_\syt(\x, \hat{\x})} = const~(w.r.t.~\x).
\end{equation*}
{A strategy we recommend} to ensure that the stationary distribution is the posterior is to decrease the strength of the regularization as the inference progresses, $\lim_{\syt \to \infty} \beta_\syt = \beta_\infty = 0$. 
The alternative {(but more difficult)} strategy is to choose the regularization potential such as to satisfy $\expectw{\hat{\x} \from \PA}{\kappa_\syt(\x, \hat{\x})}  = const~(w.r.t.~\x)$.

\begin{remark}
The regularized Fokker-Planck-Vlasov equation \cref{eq:optimal-fokker-planck-regularized} describes the evolution of the probability density for particles subject to an interacting potential $\beta_\syt\, \kappa_\syt(\x, \hat{\x})$ \cite{Duong_2023_Vlasov-interacting,Guillin_2021_Vlasov}.  The regularized drift \cref{eq:optimal-drift-regularized} gives the corresponding  McKean-Vlasov-It\^{o} process~\cref{eq:ito-process-difussion} that accounts for  particle interactions.
If the interaction potential has the form $\kappa(\x,\widehat{\x}) = V(\x-\widehat{\x})$ the potential term in  \cref{eq:optimal-fokker-planck-regularized} is the convolution  $\beta_\syt\,\nablax V \ast q_\tau$ \cite{Duong_2023_Vlasov-interacting,Guillin_2021_Vlasov}.

The addition of a regularizer, or equivalently an interaction potential, is a qualitatively correct description of the particle dynamics \cref{eq:optimal-drift-particle-system}. 
{Here, qualitatively correct means that the distribution of the unregularized analysis and regularized analysis will be the same target posterior distribution.}
The selected potential function $\kappa$ models the particle interactions in \cref{eq:optimal-drift-particle-system}.
\end{remark}
\begin{example}[Coulomb potential]
The numerical experiments in this paper use the Coulomb potential, and the corresponding repulsive electrostatic forces:
\begin{equation}
\label{eq:electrostatic-regularization}
    \kappa(\xsyt^{[e]}, \xsyt^{[i]}) = \frac{1}{\| \xsyt^{[e]} - \xsyt^{[i]} \|_2}, \qquad \nabla_{\xsyt^{[e]}} \kappa(\xsyt^{[e]}, \xsyt^{[i]}) = -\frac{\xsyt^{[e]} - \xsyt^{[i]} }{ \| \xsyt^{[e]} - \xsyt^{[i]}  \|_2^3}.
\end{equation}
Intuitively, this can be seen as a repelling force in a system of electrons (particles in our case), where the pairwise force increases with a decrease in the distance between the said pairwise particles.
{Additionally, when $\xsyt$ consists of variables at different scales (such as velocity, temperature, and salinity), $\kappa$ must be non-dimensionalized to obtain correct regularization.
One way to achieve this non-dimensionalization is to scale the variables by the inverse of the covariance (or its square root).} 
\end{example}

\subsection{Numerical time integration of particle dynamics}
\label{subsec:ts}

Consider the numerical integration of the interacting system of \cref{eq:ito-step-system}:
\begin{equation}\label{eq:ito-step}
\begin{split}
    \dif \mathfrak{X}_{\syt} &= \mathfrak{F}(\syt, \mathfrak{X}_{\syt})\, \dif \syt + \mathfrak{S}(\syt, \mathfrak{X}_{\syt})\, \dif \mathfrak{W}_\syt, \\
    \mathfrak{F}(\syt, \mathfrak{X}_\syt) &\coloneqq \big[ \widehat{\mathbf{F}}(\syt, \x^{[1]}_\syt, q_\syt)\tr,  \cdots, \widehat{\mathbf{F}}(\syt, \x^{[\Nens]}_\syt, q_\syt)\tr \big]\!\,^{\mathrm{T}}, \\
    \mathfrak{S}(\syt, \mathfrak{X}_\syt) &\coloneqq \operatorname{blkdiag}_{e = 1 \dots \Nens} \big\{ \sig(\syt, \x[e]_\syt,  q_\syt) \big\},
\end{split}
\end{equation}
where the drift term $\mathfrak{F}: \Re_+ \times \Re^{\Nstate\Nens} \to \Re^{\Nstate\Nens}$ consists of the optimal regularized drifts for each particle \cref{eq:regularized-drift} with electrostatic regularization \cref{eq:electrostatic-regularization}, the diffusion term $\mathfrak{S}: \Re_+ \times  \Re^{\Nstate\Nens} \to \Re^{\Nstate\Nens \times M\Nens}$ is a block diagonal matrix with the diffusion terms for each particle on the diagonal, and $\mathfrak{W} \syt \in \Re^{M\Nens}$ is a Wiener process. 
The time integration is challenging due to stiffness~\cite{Reich_2021_FokkerPlanck} and the presence of stochastic forcing.
For this purpose, we propose {an} implicit-explicit (IMEX) partitioning of the dynamics. 
Specifically, the drift is chosen to be the stiff component evolved implicitly, and the diffusion to be the non-stiff component evolved explicitly.
We note that an IMEX approach was also considered in~\cite{Stuart_2020_gradient-EnKF} to solve the particular case of Langevin dynamics.
Since we are interested in converging to steady state, time accuracy of the integration is not important, and so, a low order scheme will suffice. 
We restrict ourselves to considering linearly implicit methods that require one Jacobian calculation and one linear solve per step.
To this end, we consider the Rosenbrock-Euler-Maruyama (REM) scheme \cite{Hu_1996_Rosenbrock-Maruyama} where the stiff partition is evolved using the Rosenbrock-Euler method and the non-stiff partition using the Euler-Maruyama method.
Particles are advanced from $\syt$ to $\syt + \Delta_\syt$ as follows:	
\begin{equation}\label{eq:rosem}
    \mathfrak{X}_{\syt + \Delta_\syt} = \mathfrak{X}_{\syt} + \Delta_\syt\,\big(\mathbf{I}_{\Nstate\Nens} - \Delta_\syt\, \nabla_{\mathfrak{X}}\mathfrak{F}(\syt, \mathfrak{X}_{\syt}) \big)^{-1}\, \mathfrak{F}(\syt, \mathfrak{X}_{\syt}) + \sqrt{\Delta_\syt}\, \mathfrak{S}(\syt, \mathfrak{X}_{\syt})\, \boldsymbol{\xi}_{\syt},
\end{equation}
where $\boldsymbol{\xi}_{\syt} \sim \!N(\mathbf{0}, \mathbf{I}_{M\Nens})$. 
\begin{thm}\label{thm:rosenbrock-euler-maruyama}
    The Rosenbrock-Euler-Maruyama time discretization \cref{eq:rosem} for SDEs has strong order $\mathcal{O}(\Delta_\syt^{\frac{1}{2}})$ {under the assumption of a Lipschitz-continuous drift $\mathfrak{F}$.}
\end{thm}
\begin{pf}
    {Firstly, the existence and uniqueness of the solution of an SDE\cite{Kloeden_2011_sdebook} such as \cref{eq:ito-step} requires the Lipschitz continuity of the drift $\mathfrak{F}$. 
    Next,} consider the Rosenbrock-Euler-Maruayama discretization given by
    \begin{equation}\label{eq:pr1}
        \mathfrak{X}_{\syt + \Delta_\syt} = \mathfrak{X}_{\syt} + \Delta_\syt\,\big(\mathbf{I} - \Delta_\syt\, \nabla_{\mathfrak{X}}\mathfrak{F}(\syt, \mathfrak{X}_{\syt}) \big)^{-1}\, \mathfrak{F}(\syt, \mathfrak{X}_{\syt}) + \sqrt{\Delta_\syt}\, \mathfrak{S}(\syt, \mathfrak{X}_{\syt})\, \boldsymbol{\xi}_{\syt}.
    \end{equation}
    and the Neumann series expansion of $\big(\mathbf{I} - \Delta_\syt\, \nabla_{\mathfrak{X}}\mathfrak{F}(\syt, \mathfrak{X}_{\syt}) \big)^{-1}$, given by
    \begin{equation}\label{eq:pr2}
        \big(\mathbf{I} - \Delta_\syt\, \nabla_{\mathfrak{X}}\mathfrak{F}(\syt, \mathfrak{X}_{\syt}) \big)^{-1} = \mathbf{I} + \Delta_\syt\, \nabla_{\mathfrak{X}}\mathfrak{F}(\syt, \mathfrak{X}_{\syt}) + \Delta_\syt^2\, (\nabla_{\mathfrak{X}}\mathfrak{F}(\syt, \mathfrak{X}_{\syt}))^2 + \!{O}(\Delta_\syt^3).
    \end{equation}
    For the Neumann series in \cref{eq:pr2} to be convergent (and valid), we need $ \| \Delta_\syt\, \nabla_{\mathfrak{X}}\mathfrak{F}(\syt, \mathfrak{X}_{\syt}) \| < 1$ where $\| \cdot \|$ can be any matrix norm.
    Since $\mathfrak{F}(\syt, \mathfrak{X}_{\syt})$ is Lipschitz continuous, there is a constant upper bound on $\| \nabla_{\mathfrak{X}} \mathfrak{F}(\syt, \mathfrak{X}_{\syt}) \| \leq \mathfrak{K}$.
    Secondly, whenever $\Delta_\syt < \frac{1}{\mathfrak{K}}$, we have $\| \Delta_\syt\, \nabla_{\mathfrak{X}}\mathfrak{F}(\syt, \mathfrak{X}_{\syt}) \| < 1$ making the series in \cref{eq:pr2} convergent. 

    Putting \cref{eq:pr2} back in \cref{eq:pr1}, we have
    \begin{equation}\label{eq:pr3}
        \mathfrak{X}_{\syt + \Delta_\syt} = \mathfrak{X}_{\syt} + \Delta_\syt\,\mathfrak{F}(\syt, \mathfrak{X}_{\syt}) + \Delta_\syt^2 \, \nabla_{\mathfrak{X}}\mathfrak{F}(\syt, \mathfrak{X}_{\syt}) \mathfrak{F}(\syt, \mathfrak{X}_{\syt}) + \!{O}(\Delta_\syt^3) + \sqrt{\Delta_\syt}\, \mathfrak{S}(\syt, \mathfrak{X}_{\syt})\, \boldsymbol{\xi}_{\syt},
    \end{equation}
    which is equivalent to the standard Euler-Maruyama scheme with additional $\!{O}(\Delta_\syt^2)$ terms. 
    The standard Euler-Maruyama scheme is known to have strong order $\mathcal{O}(\Delta_\syt^{\frac{1}{2}})$ \cite{Kloeden_2011_sdebook}. Hence, the Rosenbrock-Euler-Maruyama {scheme} also has strong order $\mathcal{O}(\Delta_\syt^{\frac{1}{2}})$.
\end{pf}

The analytical computation of the Jacobian $\nabla_{\mathfrak{X}}\mathfrak{F}$ is expensive, since each component of the function $[\mathfrak{F}]_{e}$ $=\widehat{\mathbf{F}}(\syt, \x[e]_\syt,\mathfrak{X}_\syt)$ depends on all particles via the parameterizations of the underlying probability densities. 
It is reasonable to approximate the Jacobian by a block diagonal matrix with:
\begin{equation}
    \nabla_\mathfrak{X} \mathfrak{F}(\syt, \mathfrak{X}_\syt) \approx \operatorname{blkdiag}_{e = 1 \dots \Nens} \left\{ \nablax \widehat{\mathbf{F}}(\syt, \x,\mathfrak{X}_\syt)\big|_{\x = \x[e]_\syt} \right\},
\end{equation}
which means that a linearly implicit integration is carried out for each particle separately.
Assuming the diffusion terms are non-stiff we leave out their derivatives from the approximate Jacobian.
From \cref{eq:optimal-drift}, \cref{eq:regularized-drift}, we have the approximation
\begin{equation*}
\begin{split}
    \nablax \widehat{\mathbf{F}}(\syt, \x,\mathfrak{X}_\syt)\big|_{\x = \x[e]_\syt} &\approx
    \mathcal{A}_\syt \nabla^2_{\x,\x} \log{\PA(\x[e]_\syt)} 
    +  (\mathbf{D}(\syt, \x[e]_\syt) - \mathcal{A}_\syt) \nabla^2_{\x,\x} \log{q_\syt(\x[e]_\syt)} \\
    &\quad      -\sfrac{\beta}{\Nens} \,\mathcal{A}_\syt\, \sum_{i \neq e}  \nabla^2_{\x,\x} \kappa_\syt(\x, \xsyt^{[i]})\big|_{\x = \x[e]_\syt}.
\end{split}
\end{equation*}
To avoid the above approximation altogether, we use a finite difference approximation of Jacobian-vector products, and solve the linear system~\cref{eq:rosem} using GMRES~\cite{Saad_1986_gmres}.
We refer the reader to~\cite{Sandu_2017_analytical-JacVec} for more details. 
This is not only less expensive, but potentially captures the parametric interactions not captured by the analytical derivatives above.
However, these methods can be time-consuming for large systems, and thus, while this is a step in the right direction, the numerical solution of particle dynamics remains an open problem. 

\ifreport
\begin{remark}\label{rem:ADAM}
Since this is an optimization problem at heart, methods for performing large scale stochastic optimization such as ADAM~\cite{Kingma_2014_Adam} can also be used as an alternative, where the stiff partition is evolved with ADAM.
\end{remark}
\fi

In the limit of large ensemble sizes $\Nens \to \infty$, the expected value of the state evolves deterministically, leading to the following termination condition:
\begin{equation}
    \left\| \xmean_{\syt_* + \Delta_\syt} - \xmean_{\syt} \right\| < \epsilon\,\Delta_\syt.
\end{equation}
Specifically, the time integration is stopped at some finite time $\syt_*$ when the change in the statistical mean $\xmean_\syt$ of the particles $\Xenssyt$ relative to step-size $\Delta_\syt$ is within a desired tolerance threshold $\epsilon$, indicating the attainment of a steady state.
We also refer the reader to \cite{Crouse_2020} for other ideas on time integration of particle flows.

\section{Examples of particular VFP filters and smoothers}
\label{sec:examples-vfp}

{
As noted in~\Cref{subsec:dicretization-and-parameterization}, the drift term~\cref{eq:drift-ensemble} approximation in the VFP can be completely described by the choice of parameterization of the prior and intermediate distributions (see~\cref{tab:parameterized-families}), along with the distribution of the observation errors. 
In this section, we discuss different choices of parametric distributions.}

\subsection{Gaussian assumptions}
{
The VFP(GG) approach (see ~\Cref{subsec:dicretization-and-parameterization}) uses Gaussian assumptions on the background, $\PB(\x) = \mathcal{N}(\x \mid \xmeanb,\Pb)$, and on the current distributions, $q_\syt(\x) = \mathcal{N}(\x \mid \xmean_\syt,\P_\syt)$. If the observation errors are Gaussian, $\PO(\x)= \mathcal{N}(\x \mid \Hobs(\x) - \y,\*R)$, where $\y$ is the observation value \cref{eq:observation}, $\*R$ is the observation error covariance, and $\Hobs$ is the observation operator,  the optimal drift~\cref{eq:optimal-drift-trivial} is:
\begin{equation}\label{eq:VFP-GG-Gobs}
\begin{aligned}
    \mathbf{F}(\syt, \xsyt, q_\syt) =& 
    -(\Pb)^{-1}\,(\xsyt - \xmeanb) -\*H\tr\*R^{-1}\left(\Hobs(\xsyt) - \y\right) \\
    &-  (\mathbf{D}(\syt, \xsyt) - \*I_{\Nstate}) (\P_\syt)^{-1}\left(\xsyt - \xmean_\syt\right) + \mathbf{d}(\syt, \xsyt).
\end{aligned}
\end{equation}
Here $\*H\tr$ is the adjoint of $\Hobs$. 
If the observation errors have a  Cauchy distribution (see~\cref{tab:parameterized-families}), with $\boldsymbol{\gamma}$ the vector of scale parameters in each dimension, the optimal drift~\cref{eq:optimal-drift-trivial} is: 
\begin{equation}
\begin{aligned}
    \mathbf{F}(\syt, \xsyt, q_\syt) =& 
    -(\Pb)^{-1}\left(\*\xsyt - \xmeanb\right)
    - 2 \*H\tr ((\Hobs(\x) - \y)\oslash(\boldsymbol{\gamma}^{\circ 2} + (\Hobs(\x) - \y)^{\circ 2})) \\
    &-  (\mathbf{D}(\syt, \xsyt) - \*I_{\Nstate}) (\P_\syt)^{-1}\left(\xsyt - \xmean_\syt\right) + \mathbf{d}(\syt, \xsyt),
\end{aligned}
\end{equation}
where $\circ$ is the element-wise exponent operator and $\oslash$ is the element-wise division operator. 
}
\subsection{Laplace prior assumption}

{
The VFP(LG) filter uses a Laplace assumption on the background distribution $\PB$, and a Gaussian assumption on the intermediate distribution $q_\syt$ (see~\cref{tab:parameterized-families}). For a Gaussian observation likelihood $\PO$, the optimal drift term~\cref{eq:drift-ensemble} is:
\begin{equation}
\begin{aligned}
    \mathbf{F}(\syt, \xsyt, q_\syt) =& 
    -\frac{2}{\theta^{\rm b}}\frac{\!{K}_{\nu - 1}(\theta^{\rm b})}{\!{K}_\nu(\theta^{\rm b})} (\Pb)^{-1}\,( \xsyt - \xmeanb ) -\*H\tr\*R^{-1}\left(\Hobs(\xsyt) - \y\right) \\
    &-  (\mathbf{D}(\syt, \xsyt) - \*I_{\Nstate}) (\P_\syt)^{-1}\left(\xsyt - \xmean_\syt\right) + \mathbf{d}(\syt, \xsyt), \\
    \textnormal{where}\qquad & \theta^{\rm b} = \sqrt{2} \Vert \x - \xmeanb\Vert_{(\Pb)^{-1}}, \quad \nu = 1 - \Nstate/2.
\end{aligned}
\end{equation}
}

\subsection{Variational particle smoothing}

{
Following~\cref{eq:smoothing-analysis-gradient-loglikelihood}, we compute the drift for the strong constraint VFPS(GG) which assumes a Gaussian prior $\PB$ and a Gaussian intermediate distribution $q_\syt(\x)$. 
Under the assumption that the  observation likelihoods are Gaussian, $\PO_k(\x)= \mathcal{N}(\x \mid \Hobs_k(\x) - \y_k,\*R_k)$, the optimal drift~\cref{eq:drift-ensemble} is: 
\begin{equation*}
\begin{split}
 \mathbf{F}(\syt, \x[e]_{0,\syt}, q_\syt) &= - \mathcal{A}_\syt \left[ {(\Pb_0)}^{-1}\,( \x[e]_{0,\syt} - \xmeanb_0 )
    + \sum_{k = 1}^{K} \HH_{k}\tr\, \Mtlm_{0,k}\tr\, \mathbf{R}_{k}^{-1}\,\big( \Hobs_{k}(\Model_{0,k}(\x[e]_{0,\syt})) - \mathbf{y}_{k} \big) \right]\\
    &\quad - (\mathbf{D}(\syt, \x[e]_{0,\syt}) - \mathcal{A}_\syt) \, \*P_{0, \syt}^{-1}\left(\x[e]_{0,\syt} - \xmean_{0,\syt}\right) + \mathbf{d} (\syt, \x[e]_{0,\syt}),
\end{split}
\end{equation*}
where $\x[e]_{0,\syt}$ is one initial condition (particle).
The first two terms are the 4D-Var gradient.
The additional term $\P_{0,\syt}^{-1}\,( \x_{0,\syt} - \xmean_{0,\syt})$ pushes each particle away from the intermediate ensemble mean. 
The method VFPS(GG) is an ensemble of coupled 4D-Var runs that provides a sample of initial conditions from the posterior distribution  -- as opposed to the standard strong-constrained 4D-Var which provides only a mode of the posterior distribution.
Langevin dynamics with the Gaussian assumption on the prior -- VFPSLn(G) -- and $\Nens = 1$ leads exactly to 4D-Var preconditioned by $\mathcal{A}_\syt = \mathbf{D}(\syt, \x[e]_{0,\syt})$.
}

\ifreport
\section{Illustration of the VFP approach with regularization and diffusion}
\label{sec:illustrate-regularization-diffusion}
%
In this section, we illustrate the effects of diffusion and regularization and distribution assumptions in the VFP for a two-dimensional problem.
Specifically, we assimilate four different observations for four different background ensembles to showcase the methods (see~\cref{tab:parameterized-families}): VFP(GG), VFP(GL), VFP(LG), and VFP(LL) for different combinations of diffusion and regularization.
The background distribution, background ensemble, observation, and observation error distribution (chosen to be Gaussian) are all predefined for this example.
In \cref{fig:demo1-prior}, the small symbols represent particles and the large symbols represent observations. The shape and color of the symbols distinguish the four ensembles and their corresponding observations.

\begin{figure}[tbhp]
    \centering
    \subfloat[Prior distribution.]{
    \includegraphics[width=0.3\linewidth]{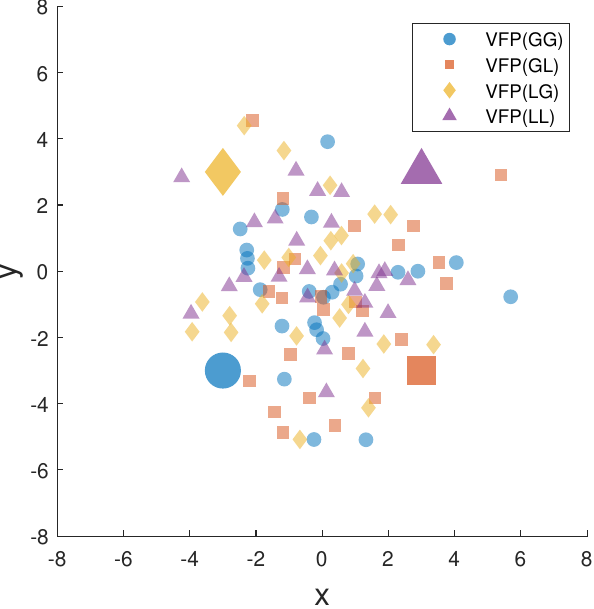}
    \label{fig:demo1-prior}
    }
    \subfloat[No diffusion or regularization.]{
    \includegraphics[width=0.3\linewidth]{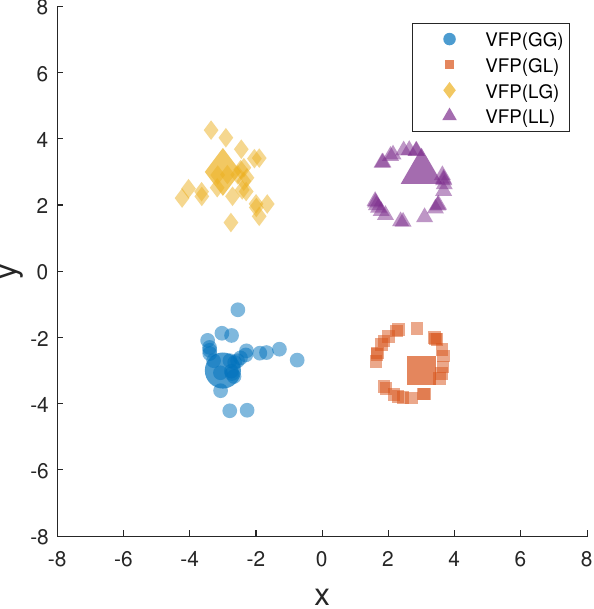}
    \label{fig:demo1-none}
    }
    \subfloat[Diffusion only.]{
    \includegraphics[width=0.3\linewidth]{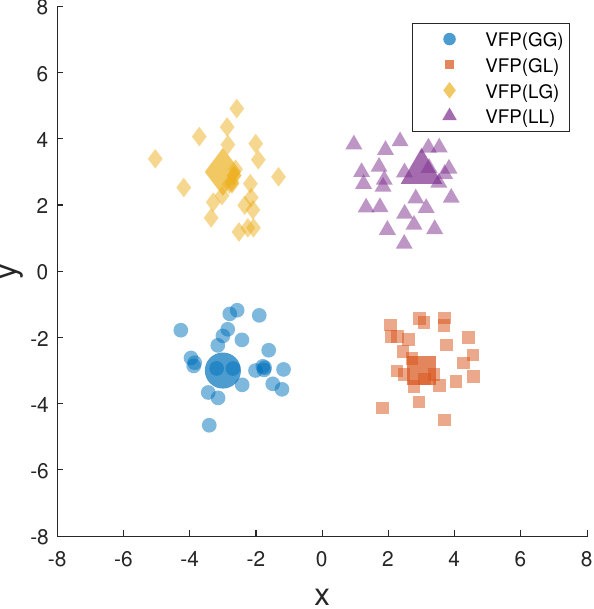}
    \label{fig:demo1-diffusion}
    }
    \\
    \subfloat[Regularization only.]{
    \includegraphics[width=0.3\linewidth]{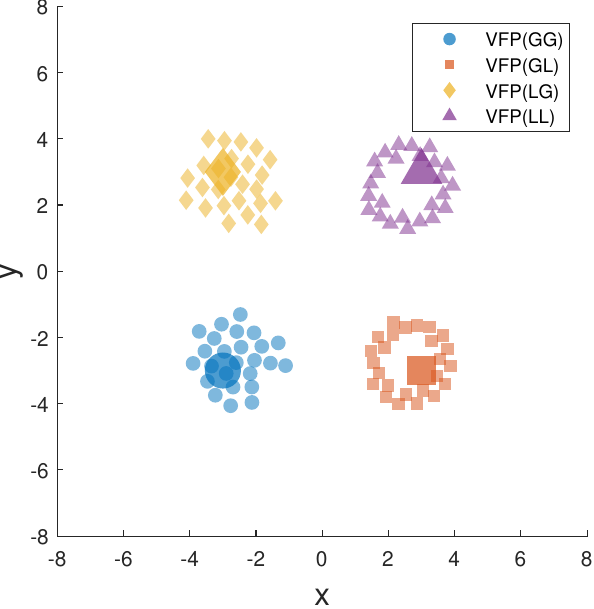}
    \label{fig:demo1-regularization}
    }
    \subfloat[Diffusion and regularization.]{
    \includegraphics[width=0.3\linewidth]{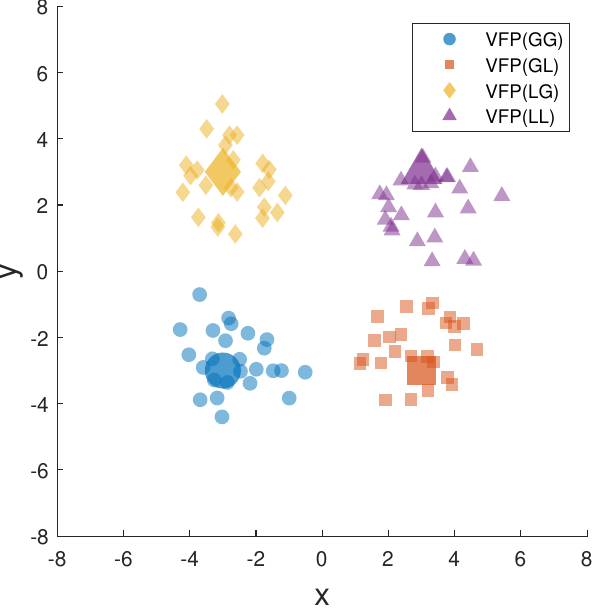}
    \label{fig:demo1-diff-and-regularization}
    }
    \caption{A comparison of the analysis distributions of VFP(GG), VFP(GL), VFP(LG), and VFP(LL) using multiple diffusion and regularization strategies. The four large symbols represent the observation corresponding to each ensemble of the same color and symbol, and the smaller symbols represent the particles.}
    \label{fig:demo-diffusion-regularization}
\end{figure}

Four experiments are performed:
\begin{enumerate*}[label={(\roman*)}]
    \item assimilation in the absence of both diffusion and regularization,
    \item assimilation in the presence of diffusion only,
    \item assimilation in the presence of regularization only,
    \item and assimilation in the presence of both diffusion and regularization.
\end{enumerate*}
and the results are shown in \cref{fig:demo-diffusion-regularization}.

\Cref{fig:demo1-none} shows the VFP analyses without diffusion or regularization. All the methods have a good number of overlapping particles, which when assimilated over multiple cycles, will result in particle collapse. 
We also observe that VFP(GL) and VFP(LL) form a ring like distribution, which is explained by the discontinuities close to the mean of the gradient of the Laplace distribution.
\Cref{fig:demo1-diffusion} shows that with the addition of diffusion, all the analysis ensembles have a diverse set of particles with no distinct pattern. 
The VFP(GG) analysis is qualitatively slightly more diverse when compared to the case without diffusion, while both VFP(LG) and VFP(LL) analyses show a much weaker ring-like structure.
With the addition of regularization in~\cref{fig:demo1-regularization}, ensembles tend to have have a much more uniform spread along their support, with each method having a unique disk-like structure. 
VFP(GL) and VFP(LL) still preserve a distinct ring-like structure while each methods shows analysis that look like non-independent sampling.
Using both diffusion and regularization in~\cref{fig:demo-diffusion-regularization}, we obtain a analyses distributions whose particles look independent and have no specific pattern. 

These results indicate that both diffusion and regularization play important roles in achieving optimal ensemble-based inference.

\fi

\section{Applying localization and covariance shrinkage in the VFP framework}
\label{sec:apply-localization}

As the prior and intermediate distributions are parameterized, certain assumptions on the distributions require the empirical estimation of the respective covariance matrices.
In high dimensional systems, when $\Nens \ll \Nstate$, the estimated covariance is of low-rank and exhibits spurious correlations between the states. 
To alleviate this effect, we use the fact that in most spatially-distributed systems the correlations between any two random states decreases with an increase in physical distance.
The flexibility of the VFP formulation allows to employ covariance  localization and shrinkage to deal with inaccurate covariance estimates. 

\subsection{Local formulation of the McKean-Vlasov-It\^{o} process}
\label{subsec:localization}

If the errors in two states $\x_{\syt, i}$ and $\x_{\syt, j}$ are conditionally independent, then the corresponding entry in the precision (inverse covariance) matrix is zero \cite{Sandu_2018_Covariance-Cholesky}. 
This implies that the Bayesian updates of a state  $\x_{\syt, i}$ depend strongly only on a subset of other variables; we denote by $\ell_i \in \{ 1, \dots, \Nstate \}$ the set of indices of these variables. 
Let $\Xens_{\syt, \{\ell_i\}} \coloneqq \Xens_{\syt}(\ell_i, :) \in \Re^{|\ell_i| \times \Nens}$ be the local ensemble of the subset of states that influence $\x_{\syt,i}$ and $\mathfrak{X}_{\syt, \{\ell_i\}} = \operatorname{vec}(\Xens_{\syt, \{\ell_i\}})$.

To perform a local update, the McKean-Vlasov-It\^{o} process must be formulated locally. 
Specifically, one computes the drift and diffusion terms in \cref{eq:ito-step} for each variable based on probability density estimates that use only local information.
Formally, for all particles $e = 1, \dots, \Nens$, the McKean-Vlasov-It\^{o} process for each state $i = 1, \dots, \Nstate$ is: 
\begin{equation}
\label{eq:coupled-ensemble-dynamics-localized}
    \dif \x_{\syt, i}^{[e]}  = \mathbf{e}_i\tr\,\mathbf{F}(\syt, \x_{\syt, \{\ell_i\}}^{[e]},\mathfrak{X}_{\syt, \{\ell_i\}})\, \dif \syt + \mathbf{e}_i\tr\,\sig(\syt, \x_{\syt, \{\ell_i\}}^{[e]},\mathfrak{X}_{\syt, \{\ell_i\}})\,\dif \*W_\syt,
\end{equation}
where $\mathbf{e}_i$ is the $i$-th standard basis vector.
This ensures that the $i$-th state is updated using only information from a local set of states, designated by the subset of indices $\{\ell_i\}$. \Cref{eq:coupled-ensemble-dynamics-localized} performs a local update of particle states.

\subsection{Local updates with Schur-product localization}
\label{subsec:schur}

Schur-product localization~\cite{Sandu_2019_adaptive-localization,Anderson_2007_localization} performs an element-wise product of the empirical covariance $\P_\syt$ with a decorrelation matrix $\mathbf{C}$ that reduces the cross covariance between distant states ~\cite{Gaspari_1999_correlation}, to obtain a localized covariance $\P_\syt^\mathrm{loc} = \P_\syt \circ \mathbf{C}$.
Each element $\mathbf{C}_{i, j} = \rho \left( d(i, j) \slash r \right)$ decreases with the physical distance $d(i, j)$ between the corresponding state variables, where $\rho(\cdot)$ is a chosen decorrelation function whose parameter $r$ is the decorrelation radius.
In our experiments, we use a combination of both local updates along with Schur-product localization using the Gaspari-Cohn decorrelation function \cite{Gaspari_1999_correlation}  with a compact support, leading to a small set of influence variables $ \ell_i = \{ j \, : \, \rho \left( d(i, j) \slash r \right) > 0 \}$, and a sparse local decorrelation matrix $\mathbf{C}_{\syt, \{\ell_i\}}$. 
Using the corresponding local anomalies $\*A_{\syt, \{\ell_i\}}$, one computes a localized covariance $\P_{\syt, \{\ell_i\}}^\mathrm{loc} = (\*A_{\syt, \{\ell_i\}} \*A_{\syt, \{\ell_i\}}\tr) \circ \mathbf{C}_{\syt, \{\ell_i\}} $\cref{eq:ensemble-mean-and-covariance} about state $i$. 
For a Gaussian assumption on the intermediate distribution, a localized version of $\nablax  \log q_\syt$ for state $i$ is 
\begin{equation}\label{eq:gaussian-intermediate-gradient-loglikelihood-loc}
        \nabla_{\x_{i}} \log q_\syt(\x_{\syt}) = - \mathbf{e}_i\tr\,(\P_{\syt, \{\ell_i\}}^\mathrm{loc})^{-1}\,(\x_{\syt, \{\ell_i\}} - \xmean_{\syt, \{\ell_i\}}).
\end{equation}
The diffusion is chosen as a scaled square root of a localized climatological covariance $\sig(\syt, \x_{\syt, \{\ell_i\}}^{[e]},\mathfrak{X}_{\syt, \{\ell_i\}}) = \diffu\, (\mathbf{B}^\mathrm{loc}_{\{\ell_i\}})^{1/2}$.
Similarly, using only the local set of states, we estimate the regularization term.
The Wiener process is different for each state $i$ and potentially lives in a smaller dimension.
This allows the states to be updated in parallel for both the filtering and smoothing cases.
We name localized VFP($\cdot, \cdot$) filters as LVFP($\cdot, \cdot$) and localized smoothers as LVFPS($\cdot, \cdot$).
However, the localization approach can quickly turn expensive due to the sheer number of computations involved. We propose an alternate method using covariance shrinkage for the Gaussian formulation.

\subsection{Covariance shrinkage}
\label{subsec:cshrink}

When the analysis and intermediate distributions in the drift~\cref{eq:drift-ensemble} are assumed Gaussian (see \cref{tab:parameterized-families}), covariance shrinkage~\cite{Chen_2010_shrinkage} can be used to alleviate the effect of spurious covariance estimates~\cite{Sandu_2015_covarianceShrinkage}. 
We have specifically used the Rao-Blackwell Ledoit-Wolf(RBLW) shrinkage estimator~\cite{Chen_2010_shrinkage} for our formulation. 
For an ensemble $\Xenssyt$, the shrinkage covariance estimate is:
\begin{align*}
    &\P_\syt^\mathrm{RBLW} = (1 - \gamma) \P_\syt + \gamma \boldsymbol{\hat{\Sigma}}_\syt, \quad
    \gamma = \min \left( \sfrac{ \left[ \frac{\Nens - 2}{\Nens} \operatorname{trace}(\P^2_\syt) \right] + \operatorname{trace}^2(\P_\syt)}{(\Nens + 2)\left[ \operatorname{trace}(\P^2_\syt) - \frac{\operatorname{trace}^2(\P_\syt)}{\Nstate}\right]} , 1\right),
\end{align*}
where the $\P_\syt $ is the covariance of intermediate particles, and
$\boldsymbol{\hat{\Sigma}}_\syt = \mu \*I_{\Nstate}$ with $\mu = \operatorname{trace}(\P_\syt)/\Nstate$ is a chosen target covariance(here, the trace-normalized identity). 
Using the Sherman-Morrison-Woodbury identity to invert $\P_\syt^\mathrm{RBLW}$, one can rewrite $\nablax  \log q_\syt$ with covariance shrinkage as: 
\begin{equation}\label{eq:gaussian-intermediate-gradient-loglikelihood-shrink}
\begin{split}
    &\nablax  \log q_\syt(\xsyt) = -(\P_\syt^\mathrm{RBLW})^{-1}\,(\xsyt - \xmean_\syt)\\
    &\qquad =  -\sfrac{(\xsyt - \xmean_\syt)}{\gamma \mu} - \sfrac{1 - \gamma}{(\gamma \mu)^2}\mathbf{A}_\syt\left( \*I_{\Nens} + \sfrac{1 - \gamma}{\gamma \mu} \*A\tr_\syt\mathbf{A}_\syt  \right)^{-1} \left( \*A\tr_\syt (\xsyt - \xmean_\syt) \right).
\end{split}
\end{equation}
We denote VFP methods that use shrinkage by ShrVFP($\cdot, \cdot$).

\section{Numerical experiments}
\label{sec:expts}

We illustrate the effectiveness of the VFP approach to data assimilation, and compare it to other state-of-the-art methods, with the help of several numerical experiments. 
All test problems and implementations are from the ODE Test problems suite~\cite{otp,otpsoft}.
{To assess the quality of the analysis and forecast ensemble, we consider the spatio-temporal root mean squared error (RMSE) and the KL-divergence of a scaled ensemble rank histogram and the ideal scaled ensemble rank histogram (i.e a uniform distribution).
The RMSE is given as
\begin{equation}\label{eq:analysis-RMSE}
    \mathrm{RMSE} = \sqrt{\sfrac{1}{\mathrm{K}\, \Nstate}\sum_{k = 1}^{\mathrm{K}} (\xt_k - \xmean_{k})\tr\,(\xt_k - \xmean_{k} )}
\end{equation}
where $\mathrm{K}$ is the number of assimilation steps (after, and excluding the initial spinup), and $\xmean_{k}$ is either the analysis ensemble mean $\xmeana_k$ (for analysis RMSE) or the forecast ensemble mean $\xmeanf_k$ (for forecast RMSE) at time $t_k$. 
Rank histograms\cite{Anderson_1996_diagram,Hamill_2001_rankhistogram} are used to assess the quality of an ensemble (analysis or forecast).
A "close to" uniform rank histogram is indicative of a good ensemble prediction as this means that each particle is equally likely to be closest to the truth.
We construct rank histograms of the first state variable of the ensemble(both analysis and forecast) taking the true trajectory as the reference. 
The second and the third state variables follow a similar trend and so, for the sake of brevity, their results are not reported here. 
To evaluate rank histograms, we define a metric that describes the KL-divergence between a discrete uniform distribution (of probabilities $\frac{1}{\Nens + 1}$) and the rank histogram values called the KLRH.
This is formalized as 
\begin{equation}\label{eq:KLRH}
    \mathrm{KLRH} = \frac{1}{\Nens + 1} \sum_{i = 1}^{\Nens + 1} \log \frac{\left(\frac{1}{\Nens + 1}\right)}{\rho_i}, 
\end{equation}
where $\rho_i$ is the scaled (analysis or forecast) ensemble ranks such that $\sum_{i = 1}^{\Nens + 1} \rho_i = 1$.
}

In all experiments, the observation error distributions are assumed to be known, and observation errors sampled from these distributions are added in \cref{eq:observation}. 

\subsection{The Lorenz '63 test problem}

For our first round of experiments, we use the 3-variable Lorenz\,'63 model~\cite{Lorenz_1963_L63},
\begin{equation}\label{eq:lorenz63}
    \sfrac{\dif x}{\dif t} = \sigma(y - x), \quad \sfrac{\dif y}{\dif t} = x(\rho - z) - y, \quad \sfrac{\dif z}{\dif t} = xy - \beta z,
\end{equation}
with the standard chaotic parameters $\sigma = 10$, $\rho = 28$, $\beta = \frac{8}{3}$.
Observations are assimilated every $\Delta t = 0.12$ time units, which is equivalent to an atmospheric time scale of 9 hours~\cite{Tandeo_2015_l63}.
{\Cref{eq:lorenz63}} is evolved in time using the Dormand-Prince 5(4) method~\cite{Dormand_1980_embeddedRK}, and the resulting discrete model is assumed to be exact, i.e $\errm = 0$. 
Each Lorenz\,'63 experiment is run for 55,000 assimilation steps discarding the first 5000 steps as the spinup. 
In the first setting, we observe all three variables with a Gaussian error coming sampled from $\!N(\*0, 8\,\*I_3)$.
In the second setting, we observe all three variables each with an independent Cauchy error sampled from $C(0, \gamma = 1)$. 
To ensure accuracy and robustness, we repeat each experiment with 12 different observation trajectories and average the results. 
Here, and in the following experiments, the observation trajectories are different noisy samples of one true model trajectory.

\subsubsection{Effect of diffusion and regularization on the ensemble rank histogram}

\begin{figure}[t]
    \centering
    \includegraphics[width=0.6\textwidth]{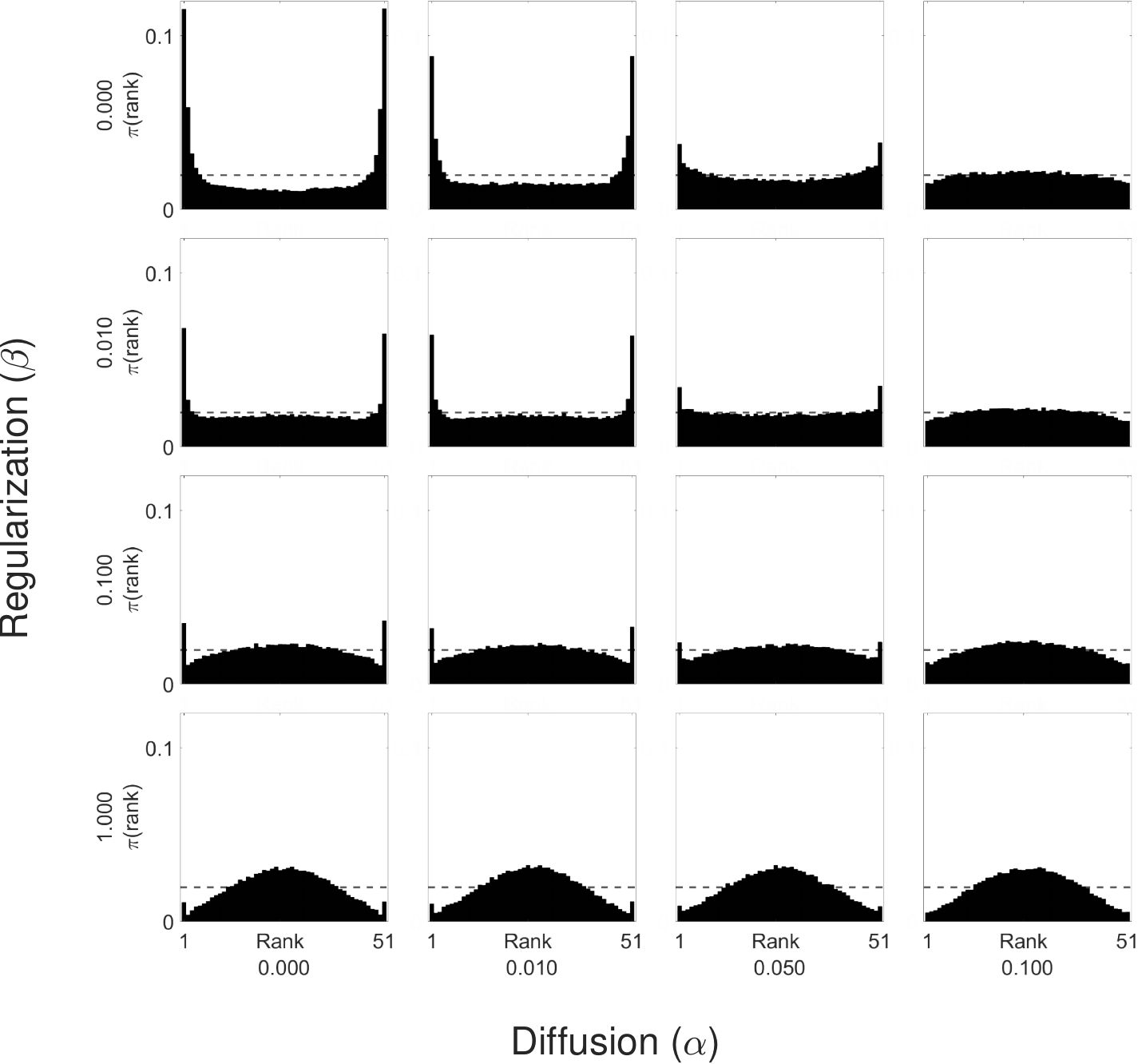}
    \caption{Lorenz-63 problem \cref{eq:lorenz63}. Analysis rank histograms for 50 particles obtained with different values of the diffusion parameter $\alpha$ and the regularization parameter $\beta$.}
    \label{fig:lorenz-63-rank-histogram}
\end{figure}

Here, we analyze the quality of the analysis ensemble for different diffusion and regularization coefficients using rank histograms\cite{Anderson_1996_diagram,Hamill_2001_rankhistogram}. 

We consider diffusion $\sig(\xsyt) = \diffu \*A^\mathrm{b}$ to be given by the background anomalies~\cref{eq:ensemble-mean-and-covariance}, and regularization by~\cref{eq:electrostatic-regularization} scaled by $\reg$. 
The assimilation scheme used in this experiment is VFP(GG), with along with the Gaussian observation error setting of $\*R = 8\,\*I_3$, and an ensemble size of $\Nens = 50$. 
The parameters for diffusion and regularization are varied as $\diffu = \{0, 0.01, 0.05, 0.1\}$ and $\reg = \{0, 0.01, 0.1, 1\}$. 
The average rank histogram of the 12 different observation trajectories for each choice of diffusion and regularization is plotted in \cref{fig:lorenz-63-rank-histogram}. 

We see is that without diffusion or regularization(leftmost row, topmost column in \cref{fig:lorenz-63-rank-histogram}), the truth often tends to fall outside the ensemble as seen in the rank histogram.
As we increase the diffusion parameter $\alpha$(i.e. along the row from left to right in \cref{fig:lorenz-63-rank-histogram}), we see that the truth starts to progressively fall more frequently inside the ensemble. 
This is explained by the diffusion increasing the spread of the ensemble. 
At the same time, if we increase the regularization parameter $\beta$(i.e. along the column from top to bottom in \cref{fig:lorenz-63-rank-histogram}), we see truth is either outside the ensemble or is pushed towards the "center" of the ensemble. 
This is explained by the choice of regularization, that acts on the ensemble.
Both the parameters $\alpha$ and $\beta$ are tuned together to obtain the optimal uniform rank histogram. 
Based on \cref{fig:lorenz-63-rank-histogram}, we choose $\alpha = 0.1$ and $\beta = 0.01$ for all subsequent experiments on the Lorenz\,'63 system, even if these values may be optimal only for this particular setting.

\subsubsection{Comparison of different VFP filters with traditional data assimilation methods}

We now compare different formulations of the variational Fokker-Planck filters corresponding to different assumptions about the underlying background and intermediate distributions. 
We consider four different formulations of the variational Fokker-Planck, namely, the VFP(KK), VFP(GG), VFP(GH), and VFPLang(G) schemes (see \cref{tab:parameterized-families} for more details about the distributions underlying these assumptions).
The Huber distribution parameters are set to $\delta_1 = \delta_2 = 1$, which were empirically found to provide good performance.
{The VFP(KK) scheme assumes a Gaussian kernel around each particle, with the kernel covariance being a scaled multivariate Scott's rule of thumb~\cite{Scott_2015_KDE} estimate given as 
\begin{equation}
    \alpha_{\rm bw} \cdot \Nens^{\frac{-2}{\Nstate + 4}} \cdot \operatorname{diag}\left( \operatorname{var}(\Xens)\right)    
\end{equation}
where $\alpha_{\rm bw}$ is hand tuned parameter, and $\operatorname{var}(\Xens)$ returns the ensemble variance of each state. 
When we tuned $\alpha_{\rm bw}$, we found that a smaller $\Nens$ required a larger $\alpha_{\rm bw}$, leading to a larger kernel covariance to compensate for the scarcity of particles.
}
We compare the performance of VFP filters against the ETKF~\cite{Bishop_2001_ETKF} (with an optimal inflation for each ensemble size $\Nens$), ETPF~\cite{Reich_2013_ETPF} and ETPF2~\cite{Reich_2017_ETPF_SOA}. 
The baseline (lower-bound) method used is the sequential importance resampling particle filter~\cite{Doucet_2001_introSMC,vanLeeuwen_2009_PF-review} (SIR) with an ensemble size of $\Nens = 1000$ particles, as in the limit $\Nens\to\infty$ it performs exact Bayesian inference.

\paragraph{Gaussian observation errors}

\begin{figure}[tbhp]
    \centering
    \subfloat[Analysis RMSE]{
    \includegraphics[width=0.43\textwidth]{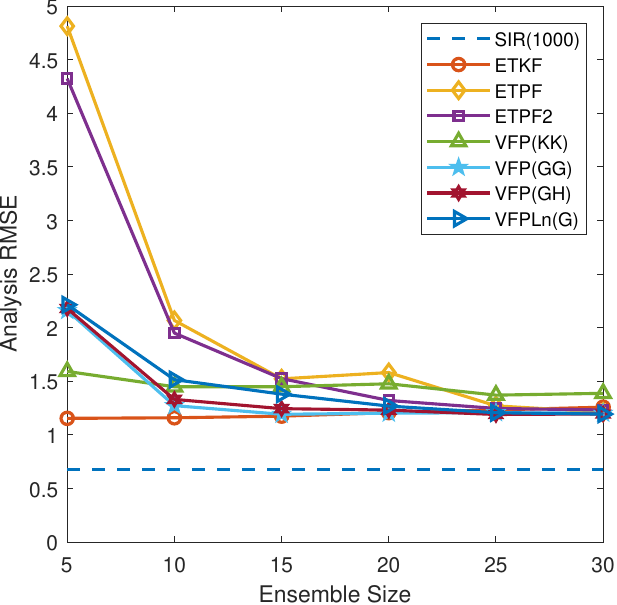}
    \label{fig:lorenz-63-gaussian-armse}
    }
    \subfloat[Forecast RMSE]{
    \includegraphics[width=0.43\textwidth]{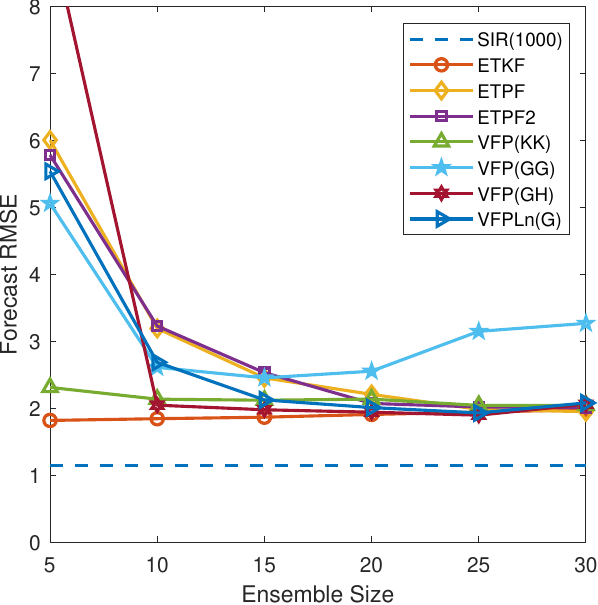}
    \label{fig:lorenz-63-gaussian-frmse}
    }
    \caption{{Lorenz-63 problem \cref{eq:lorenz63}. A comparison of analysis and forecast RMSE for multiple VFP methods along with the ETKF, ETPF, and a baseline SIR($\Nens = 1000$) with Gaussian observation error $\mathbf{R} = 8\,\*I_3$.}}
    \label{fig:lorenz-63-gaussian-results}
\end{figure}

\begin{figure}[tbhp]
    \centering
    \subfloat[Analysis RH]{
    \includegraphics[width=0.45\textwidth]{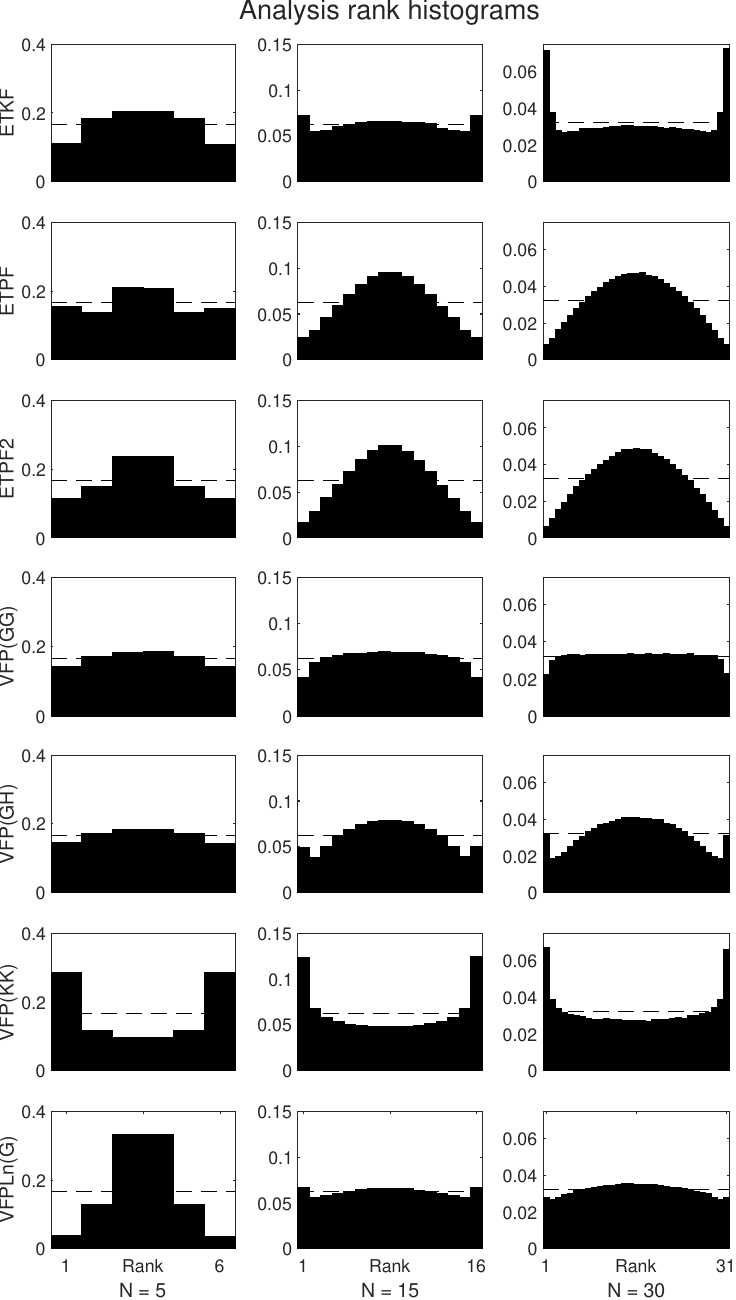}
    \label{fig:lorenz-63-gaussian-arh}
    }
    \subfloat[Forecast RH]{
    \includegraphics[width=0.45\textwidth]{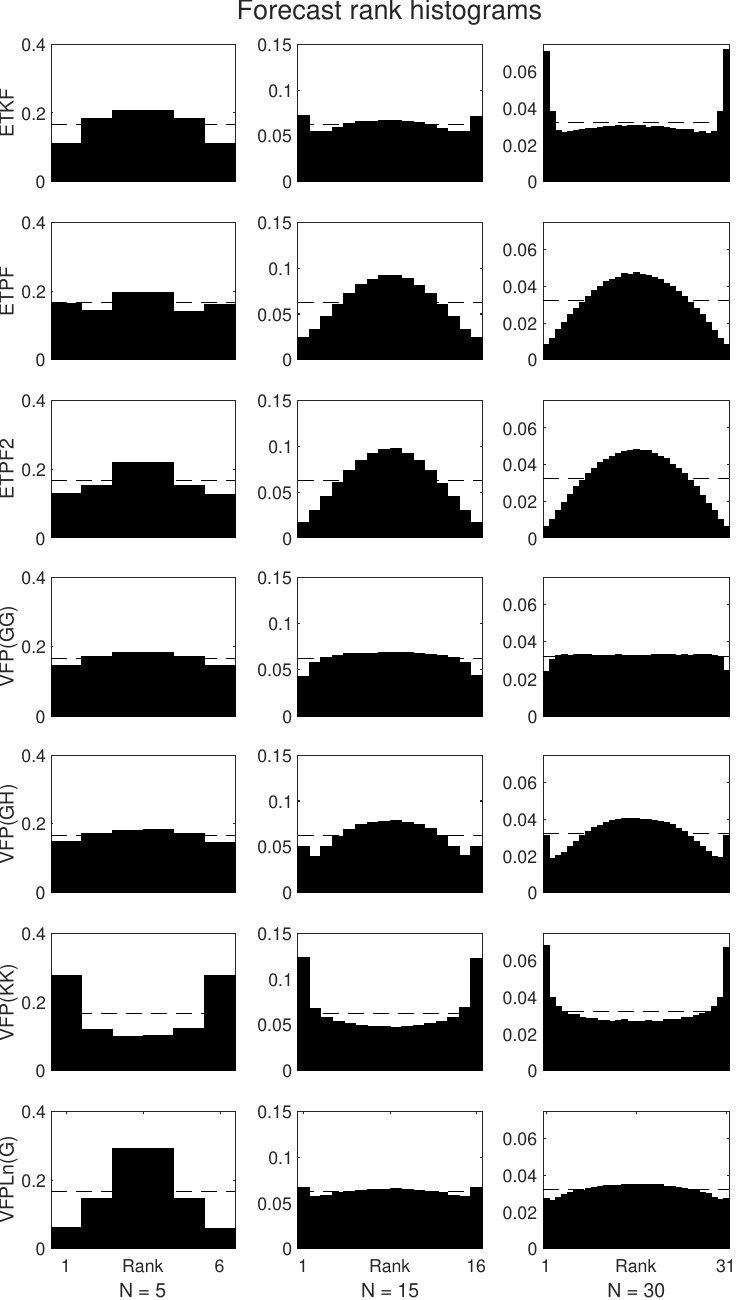}
    \label{fig:lorenz-63-gaussian-frh}
    }
    \caption{{Lorenz-63 problem \cref{eq:lorenz63}. A comparison of analysis and forecast rank histograms for multiple VFP methods along with the ETKF, ETPF with Gaussian observation error $\mathbf{R} = 8\,\*I_3$.}}
    \label{fig:lorenz-63-gaussian-results-rh}
\end{figure}

\begin{figure}[tbhp]
    \centering
    \subfloat[Analysis KLRH]{
    \includegraphics[width=0.43\textwidth]{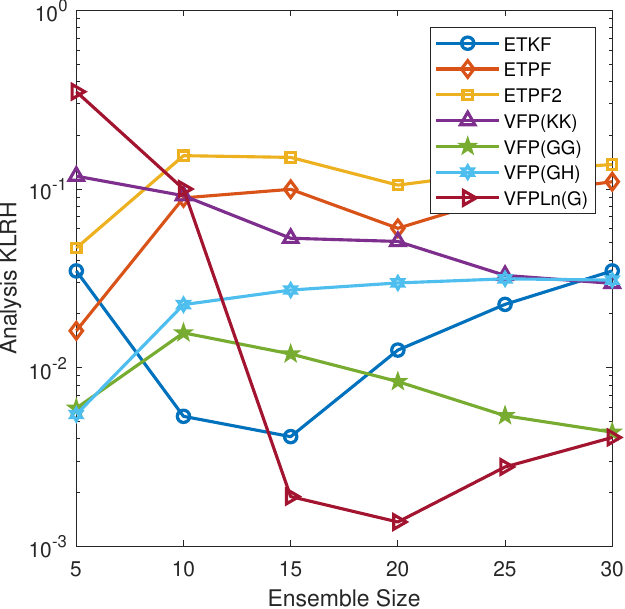}
    \label{fig:lorenz-63-gaussian-akl}
    }
    \subfloat[Forecast KLRH]{
    \includegraphics[width=0.43\textwidth]{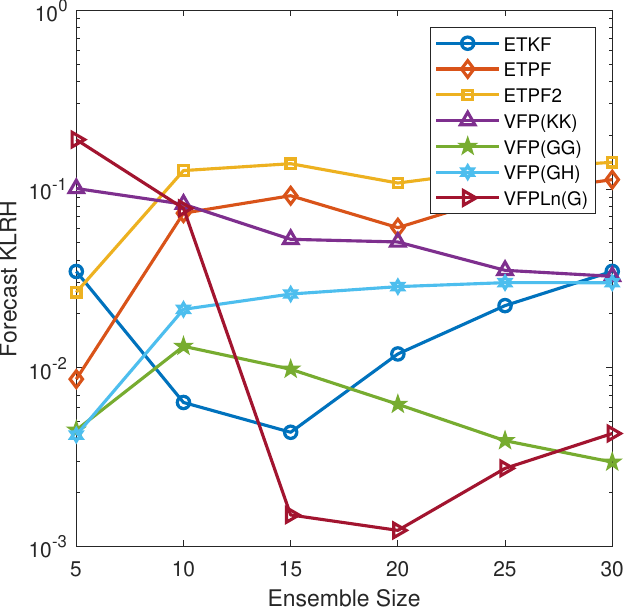}
    \label{fig:lorenz-63-gaussian-fkl}
    }
    \caption{{Lorenz-63 problem \cref{eq:lorenz63}. A comparison of analysis and forecast KLRH for multiple VFP methods along with the ETKF, ETPF with Gaussian observation error $\mathbf{R} = 8\,\*I_3$.}}
    \label{fig:lorenz-63-gaussian-results-klrh}
\end{figure}

In this setup, the observation operator is defined as $\Hobs(\x) = \x$ and unbiased Gaussian observation error covariance $\*R = 8\,\*I_3$. 
We report the analysis RMSE~\cref{eq:analysis-RMSE} for different ensemble sizes $\Nens$ in \cref{fig:lorenz-63-gaussian-armse}.
{
Both VFP(GG) and ETKF approximate fully Gaussian inference, and thus, it is not surprising that their performance is highly similar.
At lower $\Nens$ such as 5 and 10, ETKF is slightly better than VFP(GG), perhaps, purely because of parameter tuning. 
Despite the different assumption on the drift, the VFP methods show similar performance, again, due to the tuning. 
We also look at the forecast RMSE in \cref{fig:lorenz-63-gaussian-frmse}.
All methods, with the exception of VFP(GG) and VFP(GH) follow a similar trend as the analysis RMSE, albeit with a higher RMSE.
It seems as though VFP(GG) and VFP(GH) are more unstable with respect to the forecast RMSE. 
Across both the RMSEs, it seems like VFP(KK) and VFPLn(G) are the winners. 
Next, we look the KLRH (as in \cref{eq:KLRH}) for both the analysis and forecast ensembles in \cref{fig:lorenz-63-gaussian-results-klrh}.
A subset ($\Nens = 5, 15, 30$ only) of the rank histograms are depicted in \cref{fig:lorenz-63-gaussian-results-rh} to link the rank histrograms to the KLRH. 
In \cref{fig:lorenz-63-gaussian-results-rh}, we see that the analysis and forecast rank histograms for any method at a particular $\Nens$ look almost alike to the naked eye.
The slight difference is the increase in the rank value closer to the mean of the ensemble and a decrease in the rank value at the two extreme ends.
This is intuitively explained as the particles being moved closer to the truth, that increases the rank near the mean. 
Across the different methods, the VFP(GG) and VFPLn(G) have lowest KLRH (both analysis and forecast) as $\Nens \geq 15$.
These are followed by ETKF and then VFP(GH) and VFP(KK) and finally by ETPF and ETPF2.
Again, most of these results are highly sensitive to the tuning of hyperparameters that define inflation, rejuvenation, diffusion and regularization.
}

\paragraph{Cauchy observation errors}

\begin{figure}[tbhp]
    \centering
    \subfloat[Analysis RMSE]{
    \includegraphics[width=0.43\textwidth]{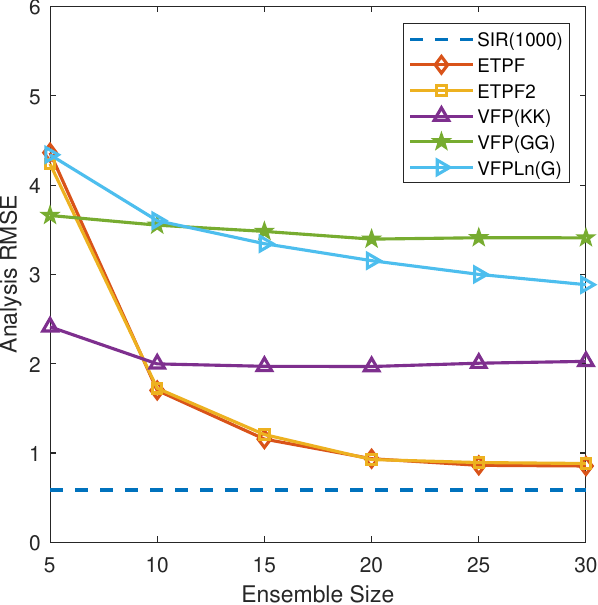}
    \label{fig:lorenz-63-cauchy-armse}
    }
    \subfloat[Forecast RMSE]{
    \includegraphics[width=0.43\textwidth]{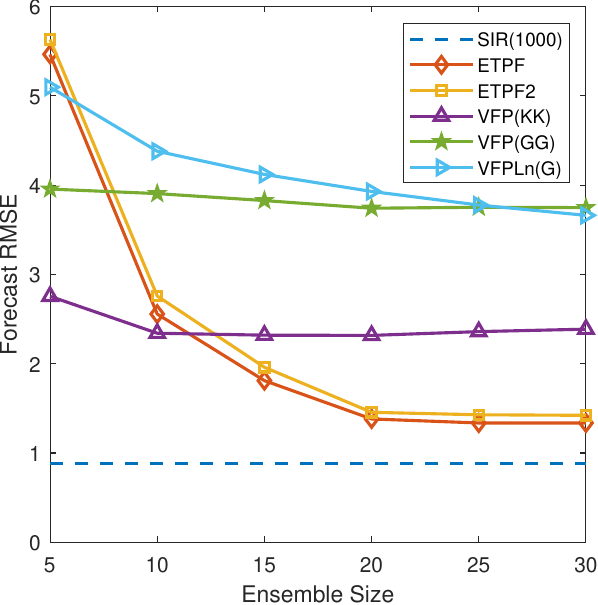}
    \label{fig:lorenz-63-cauchy-frmse}
    }
    \caption{{Lorenz-63 problem \cref{eq:lorenz63}. A comparison of analysis and forecast RMSE for multiple VFP methods along with the ETKF, ETPF, and a baseline SIR($\Nens = 1000$) with Cauchy observation error $\gamma = 1$.}}
    \label{fig:lorenz-63-cauchy-results}
\end{figure}

\begin{figure}[tbhp]
    \centering
    \subfloat[Analysis RH]{
    \includegraphics[width=0.45\textwidth]{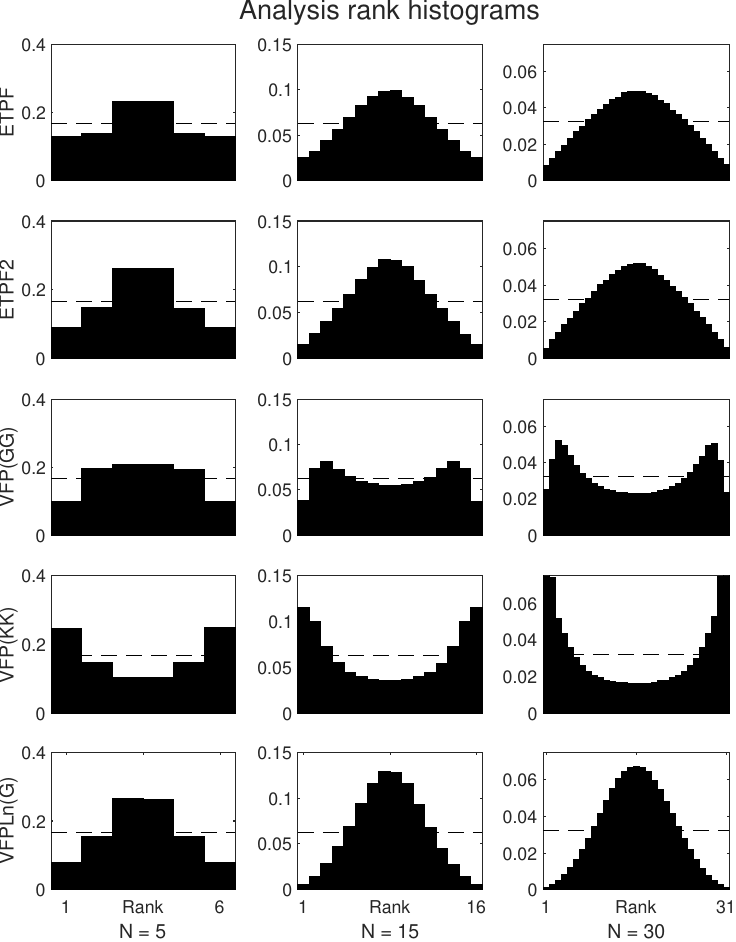}
    \label{fig:lorenz-63-cauchy-arh}
    }
    \subfloat[Forecast RH]{
    \includegraphics[width=0.45\textwidth]{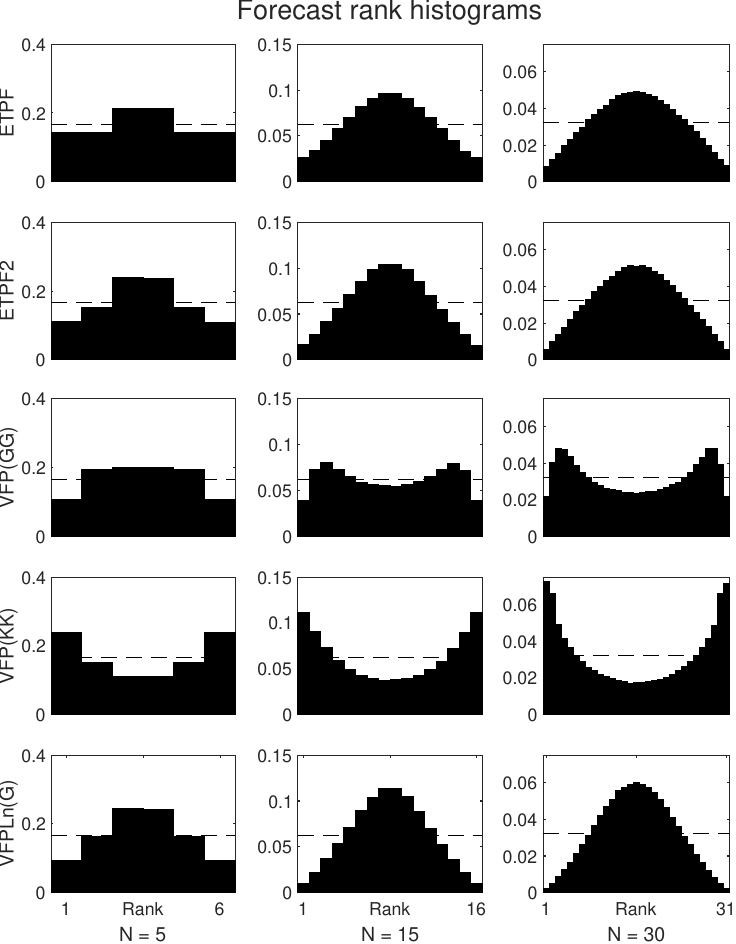}
    \label{fig:lorenz-63-cauchy-frh}
    }
    \caption{{Lorenz-63 problem \cref{eq:lorenz63}. A comparison of analysis and forecast rank histograms for multiple VFP methods along with the ETKF, ETPF with Cauchy observation error $\gamma = 1$.}}
    \label{fig:lorenz-63-cauchy-results-rh}
\end{figure}

\begin{figure}[tbhp]
    \centering
    \subfloat[Analysis KLRH]{
    \includegraphics[width=0.43\textwidth]{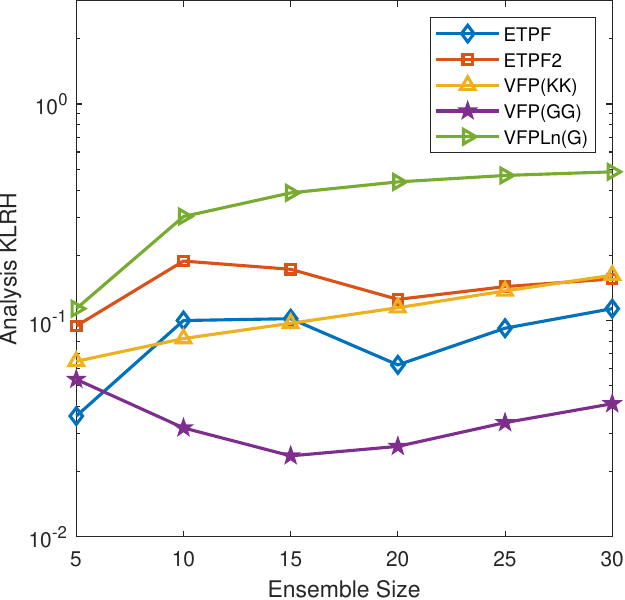}
    \label{fig:lorenz-63-cauchy-akl}
    }
    \subfloat[Forecast KLRH]{
    \includegraphics[width=0.43\textwidth]{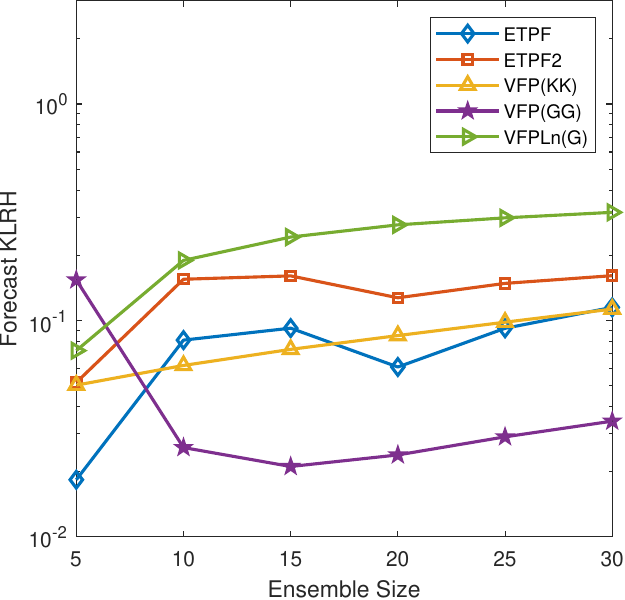}
    \label{fig:lorenz-63-cauchy-fkl}
    }
    \caption{{Lorenz-63 problem \cref{eq:lorenz63}. A comparison of analysis and forecast KLRH for multiple VFP methods along with the ETKF, ETPF with Cauchy observation error $\gamma = 1$.}}
    \label{fig:lorenz-63-cauchy-results-klrh}
\end{figure}

Here, the observation operator is, again, defined as $\Hobs(\x) = \x$ with unbiased Cauchy observation error with the parameter $\gamma = 1$ for each state.
{
The analysis RMSE for these experiments are reported in \cref{fig:lorenz-63-cauchy-armse} for various ensemble sizes $\Nens$. 
When the observations errors are sampled from the tail end of the Cauchy distribution, the flow filters end up moving the particles towards the observation that may not lie on the Lorenz '63 manifold.
Over time, this can build up and cause filter failure.
For the ETKF, {we} used with multiple choices of using an approximate $\*R$ and inflation, but the filter would ultimately diverge, and hence has not been reported in \cref{fig:lorenz-63-cauchy-armse}. 
We believe that ETKF -- based on Gaussian theory -- cannot deal with the pathological nature of the Cauchy distribution.
For the other experiments whose results are reported, some trials did fail and only the best set of 12 results were considered and reported.
VFP(GG), VFP(KK) and VFPLn(G) perform reasonably well, but are worse than the ETPF and ETPF2.
We believe that ETPF and ETPF2 show better results due to the fact that i) they make no assumptions on ensemble distributions, and ii) they optimally transport mass towards the more likely particles (essentially, enforcing bounds on how much a particle is moved towards the observation) making it robust to rare occurences of highly noisy observation errors.
VFP(KK) performs better than VFP(GG) and VFPLn(G) for the same reason, it makes a kernel density estimate of the distributions that is more accurate than a Gaussian assumption. 
The forecast RMSE in \cref{fig:lorenz-63-cauchy-frmse} demonstrate a similar trend.
Next, we look at the KLRH in \cref{fig:lorenz-63-cauchy-results-klrh} and a subset of the corresponding rank histograms in \cref{fig:lorenz-63-cauchy-results-rh}.
The rank histograms of VFP(GG) and VFP(KK) show that the methods underestimate the true state.
However, the rank histograms of ETPF, ETPF2 and VFPLn(G) overestimate the true state.  
From the forecast and analysis KLRH, it looks like VFP(GG) has the lowest values, and we believe this is mainly due to the level of diffusion.
However, all the methods have reasonably low KLRH.
}

\subsection{The Lorenz '96 test problem}

\begin{figure}[tbhp]
    \centering
    \subfloat[Lorenz '96]{
    \includegraphics[width=0.43\textwidth]{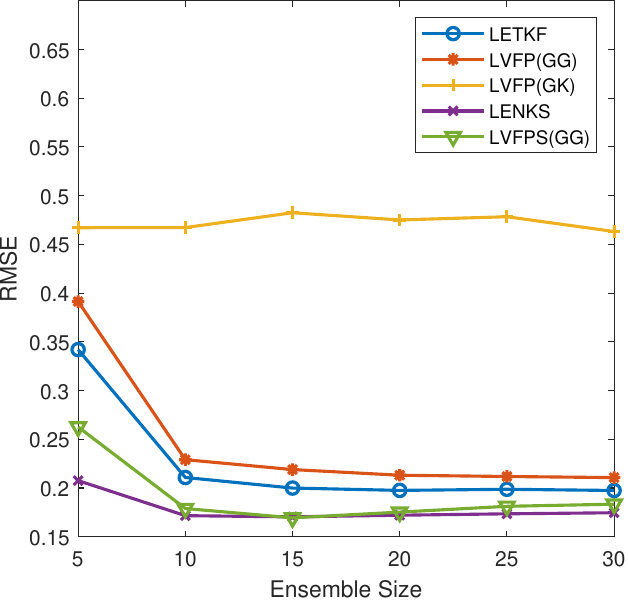}
    \label{fig:lorenz-96}
    }
    \subfloat[{Quasi-geostrophic equations}]{
    \ifreport
        \includegraphics[width=0.43\textwidth]{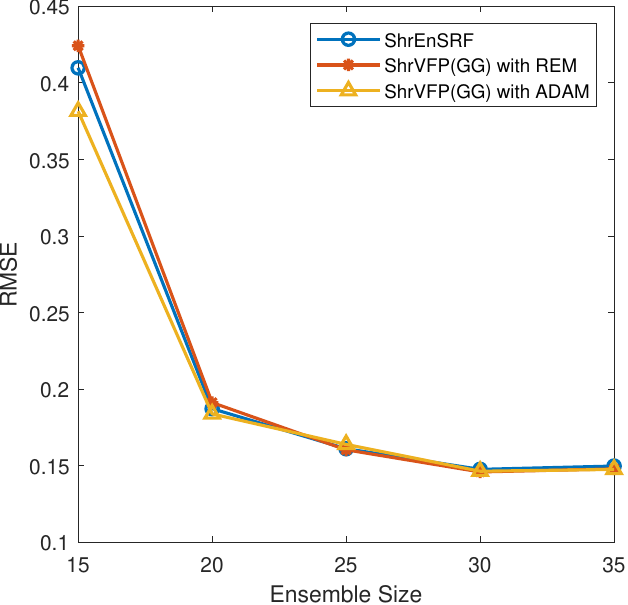} 
    \else
        \includegraphics[width=0.43\textwidth]{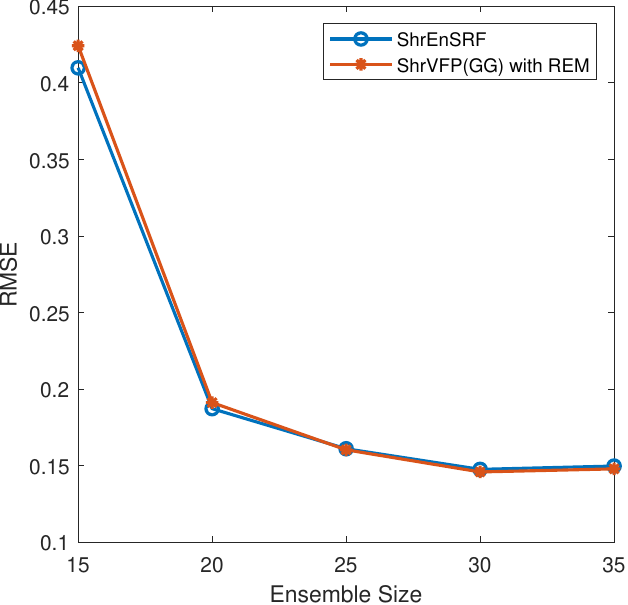}
    \fi
    \label{fig:qg}
    }
    \caption{A comparison of analysis RMSE for different ensemble sizes for the Lorenz '96 and quasi-geostrophic equations model.}
    \label{fig:lorenz96_qg}
\end{figure}

The second experiment is performed on the medium sized 40 variable Lorenz\,'96 problem~\cite{Lorenz_1996_L96,vanKekem_2018_l96dynamics} to demonstrate localized VFP(LVFP) filters and localized {strong constrained} VFP smoothers(LVFPS). 
The dynamics of the system are given by 
\begin{equation}\label{eq:L96}
    \sfrac{\dif x_i}{\dif t} = (x_{i+1} - x_{i-2})\,x_{i-1} - x_i + F, \quad \text{for} \quad i = 1, \dots, 40 \quad \text{and} \quad F = 8,
\end{equation}
where the states live on an integer ring modulo 40 i.e. -- $x_{-1} = x_{39}$, $x_0 = x_{40}$ and $x_{41} = x_1$.
We assimilate observations every $\Delta t = 0.05$ time units, which is equivalent to an atmospheric time scale of 6 hours.
The system is evolved in time using the Dormand-Prince 5(4) method~\cite{Dormand_1980_embeddedRK}, without any model error, i.e $\errm = 0$.
In the smoother, we use a discrete adjoint of the Runge-Kutta method~\cite{Sandu_2006_dadjRK}.
The experiment is run for 2200 assimilation steps where the first 200 steps are discarded as spinup. 
We observe all variables with observation operator $\Hobs(\x) = \x$.
The observation errors come from an unbiased Gaussian distribution $\!N(\*0, \*R = \*I_{40})$.
The diffusion and regularization scaling parameters are set to $\diffu = 0.1$ and $\reg = 0$.
We choose to have no regularization in these experiments to reduce the time taken to compute the optimal drift.
Each experiment is repeated with 12 different observation trajectories whose results are averaged to ensure robustness.

In \cref{fig:lorenz-96}, we compare the RMSE against $\Nens$ for three VFP methods -- namely LVFP(GG), LVFP(GK) and LVFPS(GG) -- with the localized ensemble transform Kalman filter (LETKF)~\cite{Hunt_2007_4DLETKF}, and the localized ensemble transform Kalman smoother (LETKS)~\cite{Asch_2016_book}.
In the LVFP(GG) and LVFP(GK), both Schur-localization and local update are used when evolving the intermediate ensemble.
The localization radii are fixed to be $r = 2$ for $\Nens = 5$, $r = 4$ for $\Nens = 10$ and $r = 5$ for $\Nens = 15, 20, 25, 30$ respectively in the Gaspari-Cohn decorrelation function. 
Due to the ring-like structure of Lorenz '96, $r = 2$ updates a state using information from the 7 neighboring states on either side. 
Similarly, $r = 4$ uses 14 states and $r = 5$ uses 17 states from either side.
Note that the localization we do is different from the one done in the PFF~\cite{Hu_2020_mapping-PF} manuscript.
In LVFPS(GG), only Schur-localization is performed as local updates would require a local model and local model adjoint, which we have chosen to not implement due to its impracticality for most problems of interest.
In LETKS and LVFPS(GG), we assimilate with a window size of $K = 5$ i.e. each forecast is assimilated with 5 sequential observations.
The two filters based on Gaussian assumptions(LVFP(GG) and LETKF) show a strong similarity in the RMSE. 
The same comment can be made about LVFPS(GG) and LETKS as well.
This occurs clearly because the Gaussian assumptions in LVFP(GG) and LVFPS(GG) mimic the dynamics of LETKF and LETKS that are derived from a Gaussian assumption on the ensemble.
However, LVFP(GK), had higher errors for this problem setup.
We believe that further tuning the kernels and localization, more loyal to PFF~\cite{Hu_2020_mapping-PF} could result in a better performance.
As a side note, we attempted to use other assumptions such as LVFP(GH), LVFP(HH) and LVFP(KK) for the filtering problem, whose results were unworthy to be reported. 

\subsection{The quasi-geostrophic equations test problem}

The quasi-geostrophic equations~\cite{San_2011_qg,San_2015_qge,Charney_1947_QG} approximate oceanic and atmospheric dynamics where the Coriolis and pressure gradient forces are almost balanced.
This PDE is written as: 
\begin{equation}
\label{eq:QG}
    \begin{split}
        &\mathbf{\omega}_t + \mathbf{J}(\*\psi, \*\omega) - Ro^{-1}\*\psi_x = Re^{-1} \Delta \*\omega + Ro^{-1} \mathbf{F}, \\
        &\*J(\*\psi, \*\omega) \equiv \*\psi_y \*\omega_x - \*\psi_x\*\omega_y, \quad \mathbf{F} = \sin{(\pi(y - 1))}, \quad \*\omega = - \Delta \*\psi,
    \end{split}
\end{equation}
where $\*\omega$ is the vorticity, $\*\psi$ is the streamfunction, $Re = 450$ is the Reynolds number, $Ro = 0.0036$ is the Rossby number, $\*J$ is the non-linear Jacobian, and $\mathbf{F}$ is the symmetric double gyre forcing term. 
The domain is defined to be $\Omega = [0, 1] \times [0, 2]$ which is discretized on a $63 \times 127 = 8001$ mesh. 
A constant homogeneous Dirichlet boundary condition of $\*\psi_{\partial \Omega} = 0 $ is assumed. 
{These settings result in turbulent flows with a 4 gyre circulation when averaged over time~\cite{San_2011_qg}. 
The initial condition (for the truth) is obtained by evolving a smooth random field for some time until a physically consistent field is obtained.}

\begin{figure}
    \centering
    \includegraphics[width=\linewidth]{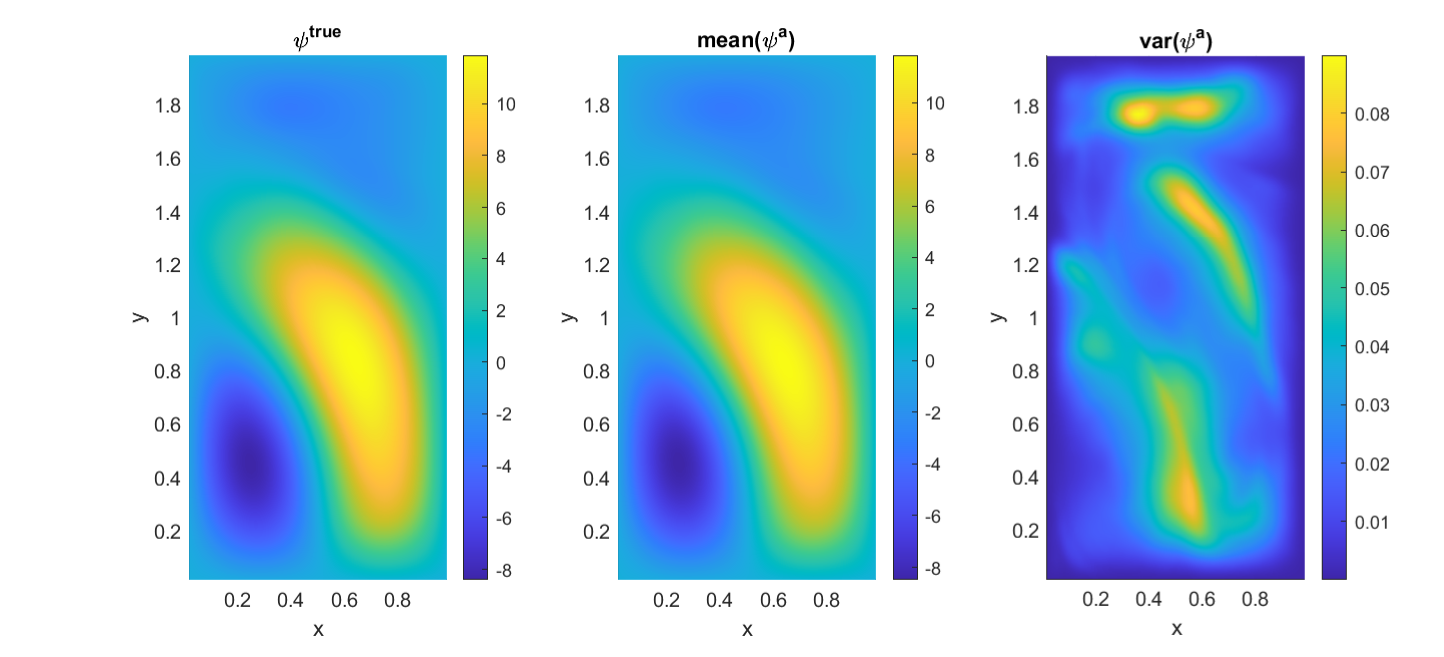}
    \caption{{An example where the figure on the left shows the true streamfunction for a simulation at 400 days. The figure in the center shows the analysis ensemble mean estimate of the streamfunction at 400 days. The figure on the right shows the ananlysis ensemble variance for the streamfunction at 400 days}}
    \label{fig:qgexample}
\end{figure}

We assimilate observations every $\Delta t = 0.0109$ time units, which is equivalent to an atmospheric time scale of 1 day.
As before, the system is evolved in time using the Dormand-Prince 5(4) method~\cite{Dormand_1980_embeddedRK}, without any model error, i.e $\errm = 0$.
The experiment is run for a total of 400 assimilation steps where the first 50 steps are discarded as spinup. 
We observe 150 evenly spaced states between 1 and 8001, which is 1.97\% of $\Nstate$.
The observation errors are sampled from a Gaussian distribution given by $\!N(\*0, \*R = \*I_{150})$.
The diffusion and regularization scaling parameters are set to $\diffu = 0.1$ and $\reg = 0$.
Again, we have no regularization in these experiments to reduce the time taken to compute the optimal drift.
As before, each experiment is repeated with 12 different observation trajectories whose results are averaged.
In \cref{fig:qg}, we compare the spatio-temporal RMSE vs the ensemble size $\Nens$.
The only method that produced a competitive result was ShrVFP(GG) which is the covariance shrinkage based VFP(GG).
LVFP(GG) was attempted, but stopped as it made very slow progress. 
To speed up the convergence, we also use the ETKF solution as the first intermediate ensemble. 
\ifreport
As time to solution was an issue, we tried evolving the system via ADAM which seemed to show competitive results as reported here.
\fi
We compare our method to a shrinkage based ensemble square root filter~\cite{Asch_2016_book} which we call ShrEnSRF. 
What we see is that VFP shows very similar performance when compared to ShrEnSRF for this problem.
This is again due to the fact that ShrEnSRF makes a Gaussian assumption similar to ShrVFP(GG). 
As stated before, further research is required to make other distribution parameterizations work well in the context of VFP in high dimensional problems. 

\section{Conclusions and future work}
\label{sec:conc}

This work discusses the Variational Fokker-Planck (VFP) framework for data assimilation, a general approach that subsumes multiple previously proposed ensemble variational methods. 
The VFP framework solves the Bayesian inference problem by smoothly transforming a set of particles into samples from the posterior distribution. 
Particles evolve in synthetic time in state-space under the flows of an ensemble of McKean-Vlasov-It\^{o} processes, and the underlying probability densities evolve according to the corresponding Fokker-Planck equations.  
We construct the optimal drift to define the McKean-Vlasov-It\^{o} processes, and show that the corresponding Fokker-Planck solutions evolve toward a unique steady-state equal to the desired posterior probability density. 
This guarantees the convergence of the VFP approach, i.e., the particles are transformed into i.i.d. samples of the posterior in the limit of infinite synthetic time. 
The choice of the diffusion terms in the McKean-Vlasov-It\^{o} processes does not change {the} evolution of the underlying probability densities toward the posterior, however it is important in practice as it acts as a particle rejuvenation approach that helps alleviate particle collapse. 

The analysis of the optimal McKean-Vlasov-It\^{o} process drift for a finite system of interacting particles leads to the conclusion that the drift contains a regularization term that nudges particles toward becoming independent random variables. 
Based on this analysis, we derive computationally-feasible approximate regularization approaches that penalize the mutual information between pairs of particles, or semi-heuristic approximations of the information.
The VFP framework can be used for both filtering and smoothing. 
We show that strong/weak-constraint VFP smoothers {(discussed in \cref{subsec:dicretization-and-parameterization-smoother})} are equivalent to ensembles of coupled strong/weak-constraint 4D-Var calculations, respectively. 
These smoothers rigorously sample the posterior distributions, unlike the popular but heuristic `ensemble of 4D-Vars' approach.

The VFP framework is very flexible and allows for implementations based on various assumptions about the type of background and intermediate distributions, e.g., Gaussian, Huber-Laplace, and Kernels. 
Moreover, localization and covariance shrinkage can be incorporated in the VFP framework to aid performance for high dimensional problems. 
We show that a semi-implicit time stepping method to solve the  McKean-Vlasov-It\^{o} processes can significantly decrease the time-to-solution in VFP, and potentially in other particle flow methods, at the expense of increased computational costs.
Numerical experiments with Lorenz '63, Lorenz '96, and the quasi-geostrophic equations test problems highlight the strengths of VFP, as well as potential areas that require further research.

Further work will investigate the choice of parameterized distributions, from the rich space of possibilities, that lead to efficient VFP implementations for high-dimensional systems.
Although the Rosenbrock-Euler-Maruyama method allows for larger timesteps, the required solution for the linear system is time consuming for high-dimensional problems. 
We will investigate more efficient time integration methods for VFP.
Localization for high dimensional problems require different drift computations for each state variable. Further investigation is required to implement localization efficiently., e.g., by taking advantage of the natural parallelization possible across state variables. 

\section*{Acknowledgement}

We thank Dr. David Higdon for helpful discussions on statistics. We thank the rest of the members of the Computational Science Laboratory at Virginia Tech for their continued help and support. 

Funding: This work was supported by the Department of Energy [ASCR DE-SC0021313]; and the National Science Foundation [CDS\&E--MSS 1953113].


\bibliographystyle{elsarticle-num}
\bibliography{Bib/references,Bib/data_assim_particle,Bib/data_assim_general,Bib/data_assim_models,Bib/data_assim_kalman,Bib/sandu}

\end{document}